\documentclass[12pt,leqno]{article}
\usepackage{amsmath}
\usepackage{amssymb}
\def\overset#1#2{{\mathrel{\mathop {{#2}_{}}\limits^{#1}}}}
\def\underset#1#2{{\mathrel{\mathop {{}_{} {#2}}\limits_{{#1}_{}}}}}
\def\upplim_#1{\underset{#1}{\overline\lim}\;}
\def\lowlim_#1{\underset{#1}{\underline\lim}\;}
\newtheorem{condition}[equation]{\indent \rm Condition}
\newtheorem{corollary}[equation]{Corollary}

\newtheorem{example}[equation]{\indent \rm {\it Example}}

\newtheorem{lemma}[equation]{Lemma}
\newtheorem{proposition}[equation]{Proposition}
\newtheorem{remark}[equation]{\indent \rm {\it Remark}}
\newtheorem{theorem}[equation]{Theorem}
\newcommand{\C}{{\mathbf{C}}}

\newcommand{\im}{{\Im\,}}
\newcommand{\Lie}{{\mathop{\mathrm{Lie}}}}

\newcommand{\NO}{NO$\frac{84}{90}$}
\renewcommand{\O}{{\mathcal{O}}}
\newcommand{\ord}{{\mathrm{ord}}}
\renewcommand{\P}{{\mathbf{P}}}
\newcommand{\Pic}{{\mathop{\mathrm{Pic}}}}
\newcommand{\Q}{{\mathbf{Q}}}
\newcommand{\R}{{\mathbf{R}}}
\newcommand{\re}{{\Re\,}}
\newcommand{\StD}{\mathrm{St}(D)}
\newcommand{\supp}{\mathrm{Supp}\,}
\newcommand{\tensor}{\otimes} 
\newcommand{\Z}{\mathbf{Z}}
\newenvironment{proof}{\par{\it Proof.}}{{\it Q.E.D.}\par\vskip3pt}
\numberwithin{equation}{section}

\title{
The Second Main Theorem for Holomorphic\\
Curves into Semi-Abelian Varieties
}
\date{\hfil}
\author{
Junjiro Noguchi, J\"org Winkelmann and Katsutoshi Yamanoi}
\begin{document}
\setlength{\baselineskip}{17pt}
\maketitle

\section{Introduction and main result.}

Let $f:\mathbf{C} \to A$ be an entire holomorphic curve
from the complex plane $\C$
into a semi-Abelian variety $A$.
It was proved by [No81] that the Zariski closure of
$f(\mathbf{C})$ is a translate of a semi-Abelian subvariety
of $A$ (logarithmic Bloch-Ochiai's theorem).
Let $D$ be an effective algebraic divisor on $A$
which is compactified to $\bar D$ on a natural
compactification $\bar A$ of $A$ (see \S3).
If $f$ omits $D$, i.e., $f(\mathbf{C}) \cap D=\emptyset$,
then $f(\mathbf{C})$ is contained in a translate of a
closed subgroup of $A$ that has no intersection with $D$
(see [No98], [SiY96]).
Note that the same holds for complex semi-tori defined
in \S3 (see [NW99]).
In particular, if $A$ is Abelian and $D$ is ample, then
$f$ is constant.
This was called Lang's conjecture.
A similar statement however is found in Bloch [Bl29],
p.\ 55, Th\'eor\`eme K without much proof, and
it is not clear what his Th\'eor\`eme K really means (cf.\ [Bl29]).

The purpose of the present paper is to establish
the quantitative version of the above result for $f$ whose
image may intersect $D$,
i.e., the second main theorem and the defect relation
(cf.\ \S\S2, 3 for the notation):
\smallskip

\noindent
{\bf Main Theorem.}
{\it
Let $f:\mathbf{C} \to M$ be a holomorphic curve
into a complex semi-torus $M$ such that the image $f(\C)$ is
Zariski dense in $M$.
Let $D$ be an effective divisor on $M$ such that
the closure $\bar D$ of $D$ in $\bar M$ is an effective divisor
on $\bar M$.
Assume that $D$ satisfies the boundary condition \ref{4.11}.
Then we have the following.
\begin{enumerate}
\item
Suppose that $f$ is of finite order $\rho_f$.
Then there is a positive integer $k_0=k_0(\rho_f, D)$
depending only on $\rho_f$ and $D$ such that
$$
T_f(r;c_1(\bar D))=N_{k_0}(r; f^*D) + O(\log r).
$$
\item
Suppose that $f$ is of infinite order.
Then there is a positive integer $k_0=k_0(f, D)$
depending on $f$ and $D$ such that
$$
T_f(r;c_1(\bar D))=N_{k_0}(r; f^*D) + O(\log T_f(r; c_1(\bar D)))
+O(\log r)||_E.
$$
\end{enumerate}
Specially, $\delta(f;\bar D)=\delta_{k_0}(f;\bar D)=0$ in both cases.
}
\medskip

See Examples \ref{4.13}, \ref{5.18}, and Proposition \ref{5.13}
that show the necessity of condition \ref{4.11}.
The most essential part of the proof is the proof of an estimate
of the proximity function (see Lemma \ref{5.1}):
\begin{equation}
\label{1.1}
m_f(r; \bar D)=O(\log r) \hbox{ or } O(\log T_f(r; c_1(\bar D)))
+O(\log r)||_E.
\end{equation}
Here, when $f$ is of infinite order, we use
one idea from R. Kobayashi [Kr98]
which is the method of the proof of Proposition (2.14) in [Kr98]
(see Lemma \ref{5.4}, (ii), and (\ref{5.10})).
The notion of logarithmic jet spaces due to [No86]
plays also a crucial role
(cf.\ [DL97] for an extension to the case of directed jets).
We then use the jet projection method developed by
[\NO], Chap.\ 6, \S3 (cf.\ [No77], [No81], and [No98]).

In \S6 we will discuss some applications of the Main Theorem.

In [Kr98] R. Kobayashi claimed (\ref{1.1}) for Abelian $A$,
but there is a part of the arguments which are heuristic,
and hard to follow rigorously.
Siu-Yeung [SiY97] claimed that for Abelian $A$
\begin{equation}
\label{1.2}
m_f(r; D) \leqq \epsilon T_f(r; c_1(D))+ O(\log r) ||_{E(\epsilon)},
\end{equation}
where $\epsilon$ is an arbitrarily given positive number,
but unfortunately there was a gap in the proof
(see Remark \ref{5.30}).
M. McQuillan [Mc96] dealt with an estimate of type (\ref{1.2}) for
some proper monoidal transformation of $\bar D \subset \bar A$ with
semi-Abelian $A$ by a method different to
those mentioned above and ours (see [Mc96], Theorem 1).

It might be appropriate at this point to recall
the higher dimensional cases in which
the second main theorem has been established.
There are actually only a few such cases that have provided
fundamental key steps.
The first was by H. Cartan [Ca33] for
$f:\mathbf{C} \to \mathbf{P}^n(\mathbf{C})$ and hyperplanes
in general position, where $\mathbf{P}^n(\mathbf{C})$ is the
$n$-dimensional complex projective space.
The Weyls-Ahlfors theory [Ah41] dealt with the same case and
the associated curves as well.
W. Stoll [St53/54] generalized the Weyls-Ahlfors theory
to the case of $f:\mathbf{C}^m \to \mathbf{P}^n(\mathbf{C})$.
Griffiths et al., [CG72], [GK73], established the second
main theorem for $f:W \to V$ with a complex affine algebraic
variety $W$ and a general
complex projective manifold $V$ such that $\mathrm{rank}\, df=\dim V$,
which was developed well by many others.
For $f:\mathbf{C} \to V$ in general, only an inequality
of the second main theorem type such as (\ref{5.29})
was proved ([No77$\sim$96], [AN91]).
Eremenko and Sodin [ES92] proved a weak second main theorem
for $f:\mathbf{C} \to \mathbf{P}^n(\mathbf{C})$ and hypersurfaces
in general position,
where the counting functions are not truncated.
In this sense, the Main Theorem adds a new case
in which an explicit second main theorem is established.

{\it Acknowledgement.}
The authors are grateful to Professor Ryoichi Kobayashi
for interesting discussions and his idea mentioned above.

\section{Order functions.}

{\bf (a)}
For a general reference of items presented in this section,
cf., e.g., [\NO].
First we recall some standard notation.
Let $\phi$ and $\psi$ be functions in a variable $r>0$
such that  $\psi >0$.
Let $E$ be a measurable subset of real positive numbers
with finite measure.
Then the expression
$$
\phi(r)= O(\psi(r))\qquad (\hbox{resp. } \phi(r)= O(\psi(r)) ||_E)
$$
stands for
$$
\upplim_{r\to\infty}\frac{|\phi(r)|}{\psi(r)} <\infty\qquad
\Big(
\hbox{resp. }
\upplim_{r\to\infty,\, r \not\in E}\frac{|\phi(r)|}{\psi(r)} <
\infty
\Big).
$$
Specially, $O(1)$ denotes a bounded term.

We use the superscript ${}^+$ to denote the positive part,
e.g., $\log^+r=\max\{0,\log r\}$.
We write $\R^+$ for the set of all real positive numbers.
We denote by $\re z$ (resp.\ $\im z$) the real (resp.\ imaginary)
part of a complex number $z\in \C$.

{\bf (b)}
Let $X$ be a compact K\"ahler manifold and let
$\omega$ be a real $(1,1)$-form on $X$.
For an entire holomorphic curve $f:\mathbf{C} \to X$
we first define the {\it order} function of $f$ with respect to $\omega$
by
$$
T_f(r;\omega)=\int_0^r\frac{dt}{t}\int_{\Delta(t)}f^*\omega ,
$$
where $\Delta(t)=\{z \in \mathbf{C};|z| < t\}$ is the
disk of radius $t$ with center at the origin of the
complex plane $\mathbf{C}$.
Let $[\omega] \in H^2(X, \mathbf{R})$ be a second cohomology
class represented by a closed real (1,1)-form $\omega$ on $X$.
Then we set
$$
T_f(r;[\omega])=T_f(r;\omega).
$$
Let $[\omega']=[\omega]$ be another representation of the class.
Since $X$ is compact K\"ahler, there is a smooth function $b$ on $X$
such that $(i/2\pi) \partial\bar\partial b=\omega' - \omega$.
There is a positive constant $C$ with $|b|\leqq C$.
Then by Jensen's formula (cf.\ [\NO], Lemma (3.39) and Remark (5.2.21))
we have
$$
|T_f(r;\omega') - T_f(r;\omega)| \leqq C.
$$
Therefore, the order function $T_f(r;[\omega])$ of $f$
with respect to the cohomology class $[\omega]$ is well-defined
up to a bounded term.
Taking a positive definite form $\omega$ on $X$,
we define the {\it order} of $f$ by
$$
\rho_f=\upplim_{r\to \infty}\;
\frac{\log T_f(r;\omega)}{\log r} \leqq \infty,
$$
which is independent of the choice of such $\omega$.
We say that $f$ is of finite order if $\rho_f < \infty$.

Let $D$ be an effective divisor on $X$.
We denote by $\supp D$ the support of $D$,
but sometimes write simply $D$ for $\supp D$ if
there is no confusion.
Assume that $f(\mathbf{C}) \not\subset D$.
Let $L(D)$ be the line bundle determined by $D$ and
let $\sigma\in H^0(X, L(D))$ be a global
holomorphic section of $L(D)$ whose divisor
$(\sigma)$ is $D$.
Take a hermitian fiber metric $\| \cdot \|$ in $L(D)$
with curvature form $\omega$, normalized so that $\omega$ represents
the first Chern class $c_1(L(D))$ of $L(D)$;
$c_1(L(D))$ will be abbreviated to $c_1(D)$.
Set
\begin{align*}
T_f(r; c_1(D)) &= T_f(r;\omega) ,  \\
m_f(r; D) &= \frac{1}{2\pi} \int_0^{2\pi}\log
\frac{1}{\|\sigma(f(re^{i\theta}))\|} d\theta.
\end{align*}
It is known that if $D$ is ample, then $f$ is rational if and
only if
$$
\lowlim_{r\to\infty}
\frac{T_f(r; c_1(\bar D))}{\log r} < \infty.
$$
One sometimes writes $T_f(r; L(D))$ for $T_f(r; c_1(D))$,
but it is noted that $T_f(r; c_1(D))$ is not depending on
a specific choice of $D$ in the homology class.
We call $m_f(r; D)$ the {\it proximity} function of $f$
for $D$.
Denoting by $\ord_z f^*D$ the order of the pull-backed divisor
$f^*D$ at $z \in \C$, we set
\begin{align*}
n(t; f^*D) &= \sum_{z \in \Delta(t)} \ord_z \; f^*D, \\
n_k(t; f^*D) &= \sum_{z \in \Delta(t)}
\min\left\{k, \ord_z \; f^*D\right\}, \\
N(r; f^*D) &= \int_1^r \frac{n(t;f^*D)}{t} dt, \\
N_k(r; f^*D) &= \int_1^r \frac{n_k(t;f^*D)}{t} dt.
\end{align*}
These are called the {\it counting} functions of $f^*D$.
Then we have the {\it F.M.T.} ({\it First Main Theorem})
(cf.\ [\NO], Chap.\ Chap. V):
\begin{equation}
\label{2.1}
T_f(r; c_1(D))=N(r;f^*D)+m_f(r;D)+O(1).
\end{equation}
The quantities
\begin{align*}
\delta(f;D) & =1 - \upplim_{r\to \infty}\;
\frac{N(r;f^*D)}{T_f(r;L(D))} \in [0,1] , \\
\delta_k(f;D) &=1 - \upplim_{r\to \infty}\;
\frac{N_k(r;f^*D)}{T_f(r;L(D))} \in [0,1]
\end{align*}
are called the {\it defects} of $f$ for $D$.

{\bf (c)}
Let $F(z)$ be a meromorphic function, and let $(F)_\infty$
(resp.\ $(F)_0$) denote the polar (resp.\ zero) divisor of $F$.
Define the proximity function of $F(z)$ by
$$
m(r, F)=\frac{1}{2\pi}\int_0^{2\pi}\log^+|F(re^{i\theta})|d\theta.
$$
Nevanlinna's order function is defined by
$$
T(r, F)=m(r, F)+N(r; (F)_\infty).
$$
Cf., e.g., [\NO], Chap.\ 6 for the basic properties of $T(r, F)$.
For instance, let $T_F(r;\omega)$ be the order function
of holomorphic $F:\C \to \P^1(\C)$ with respect to the
Fubini-Study metric form $\omega$.
Then Shimizu-Ahlfors' theorem says that
$$
T_F(r;\omega) - T(r,F)=O(1).
$$
If $F\not\equiv 0$, $T(r, 1/F)=m(r, 1/F)+N(r; (F)_0)$, and
then by Nevanlinna's F.M.T. (cf.\ [Ha64], [\NO])
\begin{equation}
\label{2.2}
T(r, F)=T\left( r, \frac{1}{F} \right) + O(1).
\end{equation}
For several meromorphic functions $F_j, 1 \leqq j \leqq l$,
on $\C$ we have
\begin{align}
\label{2.3}
T\left(r, \prod_{j=1}^l F_j\right) &\leqq \sum_{j=1}^l T(r, F_j),\\
T\left(r, \sum_{j=1}^l F_j\right) &\leqq \sum_{j=1}^l T(r, F_j)+\log l,
\nonumber\\
T(r, R(F_1, \ldots, F_l)) &\leqq O\left(\sum_{j=1}^l T(r, F_j)\right)
+O(1),
\nonumber
\end{align}
where $R(F_1, \ldots, F_l)$ is a rational function
in $F_1, \ldots, F_l$ and $R(F_1(z), \ldots, F_l(z)) \not\equiv \infty$.

\begin{lemma}
\label{2.4} {\rm (cf., [\NO], Theorem (5.2.29))}
Let $X$ be a compact K\"ahler manifold,
let $L$ be a hermitian line bundle on $X$,
and let $\sigma_1,\sigma_2\in H^0(X,L)$ with $\sigma_1\not\equiv 0$.
Let $f: \C \to X$ be
a holomorphic curve such that $f(\C) \not\subset \supp (\sigma_1)$.
Then we have
$$
T\left(r, \frac{\sigma_2}{\sigma_1}\circ f\right)
\leqq T_f(r; c_1(L))+O(1).
$$
\end{lemma}

{\it Proof.}
It follows from the definition that
$$
N\left(r; \left(\frac{\sigma_2}{\sigma_1}\circ f\right)_\infty\right)\leqq
N(r; f^*(\sigma_1)).
$$
Moreover, we have
\begin{align*}
m\left(r, \frac{\sigma_2}{\sigma_1}\circ f\right) &=
\frac{1}{2\pi} \int_{\{|z|=r\}}
\log^+ \frac{\|\sigma_2\circ f\|}{\|\sigma_1\circ f\|}d\theta \\
&\leqq
\frac{1}{2\pi} \int_{\{|z|=r\}}
\log^+ \frac{1}{\|\sigma_1\circ f\|}d\theta +O(1).
\end{align*}
Thus the required estimate follows from these and (\ref{2.1}).
{\it Q.E.D.}
\smallskip

{\bf (d)} We begin with introducing a notation for a {\it small term}.
\smallskip

{\it Definition.}
We write $S_f(r; c_1(D))$, sometimes
$S_f(r; L(D))$, to express
a small term such that 
$$
S_f(r; c_1(D))=O(\log r),
$$
if $T_f(r; c_1(D))$ is of finite order,
and
$$
S_f(r; c_1(D))=O(\log T_f(r; c_1(D)))+O(\log r)||_E,
$$
otherwise.
We use the notation $S_f(r; \omega)$ in the same sense as above
with respect to $T_f(r; \omega)$.
For a meromorphic function $F$ on $\C$, the notation
$S(r,F)$ is used to express a small term with respect to
$T(r, F)$ as well.

\begin{lemma}
\label{2.5}
{\rm (i)}  Let $F$ be a meromorphic function and let
$F^{(k)}(z)$ be the $k$-th derivative of $F$ for $k=1,2,\ldots$.
Then
$$
m\left( r, \frac{F^{(k)}}{F}\right)=S(r,F).
$$
Moreover, if $F$ is entire,
$$
T(r, F^{(k)})=T(r, F)+S(r, F), \quad k\geqq 1.
$$

{\rm (ii)} Let the notation be as in Lemma \ref{2.4}, and set
$\varphi(z)=\frac{\sigma_2}{\sigma_1}\circ f(z)$.
Suppose that $\varphi\not\equiv 0$.
Then
$$
m\left( r, \frac{\varphi^{(k)}}{\varphi}\right)=S_f(r, c_1(D)),
\quad k \geqq 1.
$$
\end{lemma}

{\it Proof.}
The item (i) is called Nevanlinna's lemma on logarithmic
derivatives (cf.\ [\NO], Corollary (6.1.19)).
Then (ii) follows from (i) and Lemma \ref{2.4}.
{\it Q.E.D.}
\smallskip

The following is called Borel's lemma (cf.\ [Ha64], p.\ 38, Lemma 2.4).
\begin{lemma}
\label{2.6}
Let $\phi(r)$ be a continuous, increasing function on $\R^+$
such that $\phi(r_0)>0$ for some $r_0\in\R^+$.
Then we have
$$
\phi\left(r+\frac{1}{\phi(r)}\right)<2\phi(r) ||_E.
$$
\end{lemma}

For a later use we show

\begin{lemma}
\label{2.7}
Let $F$ be an entire function, and let $0<r<R$.
\begin{enumerate}
\item
$\displaystyle T(r, F)=m(r, F)\leqq \max_{|z|=r}\log |F(z)|
\leqq \frac{R+r}{R-r}m(R, F)$.
\item
$m(r, F)=S(r, e^{2\pi i F})$.
\end{enumerate}
\end{lemma}

{\it Proof.}
(i) See [\NO], Theorem (5.3.13) or [Ha64], p.\ 18, Theorem 1.6.

(ii) Using the complex Poisson kernel, we have
$$
F(z)=\frac{i}{2\pi}\int_{\{|\zeta|=R\}}
\frac{\zeta+z}{\zeta-z}\im F(\zeta)d\theta
+\re F(0).
$$
Therefore, using (i) and the F.M.T. (\ref{2.2}) with
$0 < r < R < R'$ we obtain
\begin{align}
\label{2.8}
\max_{|z|=r}|F(z)|&
\leqq \frac{R+r}{R-r}\max_{|\zeta|=R}|\im F(\zeta)|+|\re F(0)|\\
&\leqq \frac{R+r}{R-r}\left(
\max_{|\zeta|=R} \im F(\zeta) +
\max_{|\zeta|=R} \im(-F(\zeta))\right)
+|\re F(0)|
\nonumber\\
&\leqq \frac{R+r}{R-r}\cdot
\frac{R'+R}{R'-R}\cdot\frac{1}{2\pi}
\Big(m(R', e^{-2\pi i F})+m(R', e^{2\pi i F})
\Big)+|\re F(0)|
\nonumber\\
&\leqq \frac{R+r}{R-r}\cdot
\frac{R'+R}{R'-R}\cdot\frac{1}{\pi}
\left( T(R', e^{2\pi i F}) +O(1)\right)+O(1).
\nonumber
\end{align}
If $T(r, e^{2\pi i F})$ has finite order,
then setting $R=2r$ and $R'=3r$, we see by (\ref{2.8}) that
$$
m(r, F)\leqq \log \max_{|z|=r} |F(z)| =O(\log r).
$$
In the case where $T(r, e^{2\pi i F})$ has infinite order,
we write $T(r)=T(r, e^{2\pi i F})$ for the sake of simplicity.
Setting $R=r+\frac{1}{2T(r)}$
and $R'=r+\frac{1}{T(r)}$,
we have by (\ref{2.8}) and Lemma \ref{2.6}
\begin{align*}
m(r, F) & \leqq \log \max_{|z|=r}|F(z)|\\
&\leqq \log\left( (4rT(r)+1)(4rT(r)+3)
\left(T\left(r+\frac{1}{T(r)}\right) +O(1)\right) +O(1)\right)\\
&\leqq \log\Big(
(4rT(r)+1)(4rT(r)+3)(2T(r) +O(1)) +O(1)\Big)||_E\\
&=S(r, e^{2\pi i F}).
\end{align*}
{\it Q.E.D.}

\section{Complex semi-torus.}

Let $M$ be a complex Lie group admitting the exact sequence
\begin{equation}
\label{3.1}
0 \to (\mathbf{C}^*)^p \to M\, \overset{\eta}{\to}\, M_0 \to 0,
\end{equation}
where $\mathbf{C}^*$ is the multiplicative group
of non-zero complex numbers, and $M_0$ is a (compact) complex
torus.
Such $M$ is called a {\it complex semi-torus} or a
{\it quasi-torus}.
If $M_0$ is algebraic, that is, an Abelian variety,
$M$ is called a {\it semi-Abelian variety} or
a {\it quasi-Abelian variety}.
In this section and in the next, we assume that
$M$ is a complex semi-torus.

Taking the universal coverings of (\ref{3.1}), one gets
\begin{equation}
\nonumber
0 \to \mathbf{C}^p \to \mathbf{C}^n \to \mathbf{C}^m \to 0,
\end{equation}
and an additive discrete subgroup $\Lambda$ of $\mathbf{C}^n$ such that
\begin{align}
\nonumber
\pi :  \mathbf{C}^n &\to M =\mathbf{C}^n/\Lambda,\\
\pi_0 :  \mathbf{C}^m =(\mathbf{C}^n/\mathbf{C}^p) &\to
M_0=(\mathbf{C}^n/\mathbf{C}^p)/(\Lambda/\mathbf{C}^p), \nonumber \\
(\mathbf{C}^*)^p & = \mathbf{C}^p/(\Lambda \cap \mathbf{C}^p).
\nonumber
\end{align}
We fix a linear
complex coordinate system
$x=(x',x'')=(x'_1,\ldots,x'_p,x''_1, \ldots, x''_m)$
on $\C^n$
such that $\mathbf{C}^p \cong \{x''_1=\cdots=x''_m=0\}$
and
$$
\Lambda \cap \mathbf{C}^p=\mathbf{Z}
\begin{pmatrix}
1 \\ \vdots \\ 0
\end{pmatrix}
+\cdots+
\mathbf{Z}
\begin{pmatrix}
0 \\ \vdots \\ 1
\end{pmatrix}.
$$

The covering mapping $\mathbf{C}^p \to (\mathbf{C}^*)^p$ is given by
$$
\begin{pmatrix}
x_1 \\ \vdots \\ x_p
\end{pmatrix}
\in \mathbf{C}^p
\to
\begin{pmatrix}
e^{2\pi i x_1} \\ \vdots \\ e^{2\pi i x_p}
\end{pmatrix}=
\begin{pmatrix}
u_1 \\ \vdots \\ u_p
\end{pmatrix}
\in (\mathbf{C}^*)^p.
$$

We may regard $\eta:M \to M_0$ to be a flat
$(\mathbf{C}^*)^p$-principal fiber bundle.
By a suitable change of coordinates $(x'_1,\ldots,x'_p,x''_1, \ldots,
x''_m)$
the discrete group $\Lambda$ is generated over
$\Z$ by the column vectors of the matrix of the following type
\begin{equation}
\label{3.2}
\begin{pmatrix}
1 &  \cdots  & 0 &   \\
\vdots & \ddots & \vdots &  A \\
0 & \cdots  & 1 &                \\
 & O       &    &  B
\end{pmatrix},
\end{equation}
where $A$ is a {\it real} $(p, m)$-matrix and $C$ is a $(m,m)$-matrix.
Therefore the transition matrix-functions of the flat
$(\mathbf{C}^*)^p$-principal fiber bundle
$\eta:M \to M_0$ are expressed by a diagonal matrix
such that
\begin{equation}
\label{3.3}
\begin{pmatrix}
a_1 & 0 & \cdots & 0 \\
0  & a_2 & \cdots& 0 \\
\vdots & \vdots & \ddots & \vdots \\
0 & 0  & \cdots & a_p
\end{pmatrix},
\qquad |a_1|= \cdots =|a_p|=1.
\end{equation}

Taking the natural compactification
$\mathbf{C}^*=\mathbf{P}^1(\mathbf{C})\setminus\{0, \infty\}
\hookrightarrow \mathbf{P}^1(\mathbf{C})$,
we have a compactification of $M$,
$$
\bar\eta:\bar M \to M_0,
$$
which is a flat $(\mathbf{P}^1(\mathbf{C}))^p$-fiber bundle
over $M_0$.
Set
$$
\partial M=\bar M \setminus M ,
$$
which is a divisor on $\bar M$ with only simple
normal crossings.

Let $\Omega_1$ be the product of
the Fubini-Study metric forms on
$(\mathbf{P}^1(\mathbf{C}))^p$,
$$
\Omega_1=\frac{i}{2\pi} \sum_{j=1}^p
\frac{du_j\wedge d\bar u_j}{(1+|u_j|^2)^2}.
$$
Because of (\ref{3.3}) $\Omega_1$ is defined well
on $\bar M$.
Let $\Omega_2=(i/2\pi)\partial\bar\partial\sum_j|x_j''|^2$
be the flat hermitian metric form on $\C^m$,
and as well on the complex torus $M_0$.
Then we set
\begin{equation}
\label{3.4}
\Omega=\Omega_1+\bar\eta^*\Omega_2,
\end{equation}
which is a K\"ahler form on $\bar M$.
\medskip

{\it Remark.}
The same complex Lie group $M$ may admit several such
exact sequences as (\ref{3.1})
which may be quite different.
For instance, let $\tau$ be an arbitrary complex number with
$\im \tau >0$.
Let $\Lambda$ be the discrete subgroup of $\C^2$ generated by
$$
\begin{pmatrix} 1 \\ 0 \end{pmatrix},
\begin{pmatrix} 0 \\ 1 \end{pmatrix},
\begin{pmatrix} i \\ \tau \end{pmatrix}.
$$
Then $M=\C^2/\Lambda$ is a complex semi-torus and the natural projection
of $\C^2=\C\times\C$ onto the first and the second factors
induce respectively exact sequences of the forms
\begin{align*}
0 \to \C^* \to & M \to \C/\!\left<1,i\right>_{\Z} \to 0, \\
0 \to \C^* \to & M \to \C/\!\left<1,\tau\right>_{\Z} \to 0.
\end{align*}
In the sequel we always consider a complex semi-torus $M$ with a {\em fixed}
exact sequence as in (\ref{3.1}) and with the discrete
subgroup $\Lambda$ satisfying (\ref{3.2}).
\medskip

Let $f: \mathbf{C} \to M$ be a holomorphic curve.
We regard $f$ as a holomorphic curve into $\bar M$
equipped with the K\"ahler form $\Omega$, and define
the order function by
$$
T_f(r;\Omega)=\int_0^r\frac{dt}{t}\int_{\Delta(t)}f^*\Omega .
$$
Let $\tilde f:\mathbf{C}\to \mathbf{C}^n$ be the lifting
of $f$, and set
$$
\tilde f (z)=\left(F_1(z), \ldots , F_p(z),
G_1(z), \ldots, G_m(z) \right),
$$
where $F_i(z)$ and $G_j(z)$ are entire functions.
Extending the base $M_0$ of the fiber bundle $M \to M_0$ to
the universal covering $\pi_0:\mathbf{C}^m \to M_0$, we have
\begin{equation}
\nonumber
M \times_{M_0} \mathbf{C}^m \cong (\mathbf{C}^*)^p \times \mathbf{C}^m,
\qquad
\bar M \times_{M_0} \mathbf{C}^m \cong
(\mathbf{P}^1(\mathbf{C}))^p \times \mathbf{C}^m .
\end{equation}
Set
\begin{equation}
\nonumber
\hat M=(\mathbf{P}^1(\mathbf{C}))^p \times \mathbf{C}^m.
\end{equation}
Then $\hat M$ is the universal covering of $\bar M$,
and then $\tilde f$ induces a lifting $\hat f$ of
$f:\C \to M \hookrightarrow \bar M$,
\begin{equation}
\nonumber
\hat f: z \in \mathbf{C} \to
\left(e^{2\pi i F_1(z)}, \ldots, e^{2\pi i F_p(z)},
G_1(z), \ldots, G_m(z)\right) \in
(\mathbf{C}^*)^p \times \mathbf{C}^m=\hat M.
\end{equation}
Set
\begin{align}
\nonumber
\hat f_{(1)} &:  z \in \mathbf{C} \to
\left(e^{2\pi i F_1(z)}, \ldots, e^{2\pi i F_p(z)}
\right) \in (\mathbf{C}^*)^p , \\
\hat f_{(2)} &: z \in \mathbf{C} \to
\left(G_1(z), \ldots, G_m(z)\right) \in
\mathbf{C}^m. \nonumber
\end{align}
By definition we have
\begin{equation}
\label{3.5}
T_f(r; \Omega) =
T_{\hat f_{(1)}}(r; \Omega_1)+ T_{\hat f_{(2)}}(r; \Omega_2).
\end{equation}
By Shimizu-Ahlfors' theorem we have
\begin{equation}
\label{3.6}
T_{\hat f_{(1)}}(r; \Omega_1)=\sum_{j=1}^p T(r, e^{2\pi i F_j})+O(1).
\end{equation}
By Jensen's formula (cf.\ [\NO], Lemma (3.3.17)) we have
\begin{align}
\label{3.7}
T_{\hat f_{(2)}}(r; \Omega_2) &=
\int_0^r \frac{dt}{t}
\int_{\Delta(t)} \frac{i}{2\pi} \partial \bar\partial
\sum_{j=1}^m |G_j(z)|^2 \\
&=\frac{1}{4\pi} \int_0^{2\pi}
\left(
\sum_{j=1}^m \left| G_j (re^{i\theta})\right|^2
\right) d\theta
-\frac{1}{2} \sum_{j=1}^m \left| G_j(0)\right|^2.
\nonumber
\end{align}

\begin{lemma}
\label{3.8}
Let the notation be as above.
Then for $k \geqq 0$ we have
\begin{align*}
T(r, F_j^{(k)}) &=T(r, F_j)+ kS(r, F_j)\leqq
S_{\hat f_{(1)}}(r; \Omega_1)\leqq S_f(r; \Omega),\\
T(r, G_j^{(k)}) &=T(r, G_j)+kS(r, G_j)
\leqq S_{\hat f_{(2)}}(r; \Omega_2)\leqq S_f(r; \Omega).
\end{align*}
\end{lemma}

{\it Proof.}
By Lemma \ref{2.5} it suffices to show the case of $k=0$.
By Lemma \ref{2.7} and (\ref{3.6})
$$
T(r, F_j)=m(r, F_j)=S(r, e^{2\pi i F_j})\leqq S_{\hat f_{(1)}}(r;
\Omega_1)
\leqq S_f(r; \Omega).
$$
For $G_j$ we have by making use of (\ref{3.7})
and the concavity of the logarithmic function
\begin{align*}
T(r; G_j) &=m(r, G_j)=\frac{1}{2\pi}\int_{\{|z|=r\}}
\log^+ |G_j(z)| d\theta\\
&=\frac{1}{4\pi}\int_{\{|z|=r\}}\log^+ |G_j(z)|^2 d\theta \\
&\leqq \frac{1}{4\pi}\int_{\{|z|=r\}}\log(1+ |G_j(z)|^2) d\theta \\
&\leqq \frac{1}{2}\log\left(
1+\frac{1}{2\pi}\int_{\{|z|=r\}} |G_j(z)|^2 d\theta\right) \\
&=S_{\hat f_{(2)}}(r; \Omega_2)\leqq S_f(r; \Omega).
\end{align*}
{\it Q.E.D.}
\medskip

\begin{lemma}
\label{3.9}
Let the notation be as above.
Assume that $f:\mathbf{C} \to M$ has a finite order $\rho_f$.
Then $F_j(z), 1 \leqq j \leqq p$, are polynomials of
degree at most $\rho_f$, and $G_k, 1 \leqq k \leqq m$, are
polynomials of degree at most $\rho_f/2$;
moreover, at least one of $F_j$ has degree $\rho_f$, or
at least one of $G_k$ has degree $\rho_f/2$.
\end{lemma}

{\it Proof.}
Let $\epsilon >0$ be an arbitrary positive number.
Then there is a $r_0>0$ such that
$$
T_f(r; \Omega) \leqq r^{\rho_f+\epsilon}, \qquad r \geqq r_0.
$$
It follows from (\ref{3.5})$\sim$(\ref{3.7}) that for $r \geqq r_0$
\begin{align}
\label{3.10}
T_{\hat f_{(1)}}(r; \Omega_1) &\leqq r^{\rho_f+\epsilon}, \\
\nonumber
\frac{1}{2\pi} \int_0^{2\pi}
\left(
\sum_{j=1}^m \left| G_j (re^{i\theta})\right|^2
\right) d\theta
&\leqq r^{\rho_f+\epsilon}.
\end{align}
It follows from (\ref{3.10}), (\ref{3.6}), and (\ref{2.8}) applied
with $R=2r$ and $R'=3r$ that there is a positive constant $C$
such that
$$
\max_{|z|=r} |F(z)| \leqq C r^{\rho_f+\epsilon}.
$$
Therefore, $F_j(z)$ is a polynomial of degree at most $\rho_f$.

Expand $G_j(z)=\sum_{\nu}^\infty c_{j\nu}z^\nu$.
Then one gets
$$
\frac{1}{2\pi} \int_0^{2\pi}
\left(
\sum_{j=1}^m \left| G_j (re^{i\theta})\right|^2
\right) d\theta=
\sum_{j=1}^m \sum_{\nu=0}^\infty
|c_{j\nu}|^2 r^{2\nu}.
$$
It follows that
$$
\sum_{j=1}^m \sum_{\nu=1}^\infty
|c_{j\nu}|^2 r^{2\nu} \leqq r^{\rho_f+\epsilon},
\quad r \geqq r_0.
$$
Hence, $c_{j\nu}=0$ for all $\nu > \rho_f/2$ and $1 \leqq j \leqq m$.
We see that $G_j(z), 1 \leqq j \leqq m$, are polynomials of
degree at most $\rho_f/2$.
The remaining part is clear.
{\it Q.E.D.}
\medskip

In the language of Lie group theory we obtain the following
characterization of holomorphic curves of finite order:
\begin{proposition}
\label{3.11}
Let $M$ be an $n$-dimensional complex semi-torus
with the above compactification
$\bar M$, let $\Lie(M)$ be its Lie algebra,
and let $\exp:\Lie(M)\to M$ be the exponential map.
Let $f:\C\to M$ be a holomorphic curve.
Then $f$ is of finite order considered as a holomorphic
curve into $\bar M$ if and only
if there is a polynomial map $P:\C\to\Lie(M)\cong \C^n$
such that $f=\exp\circ P$, and hence
the property of $f$ being of finite order
is independent of the choice of the compactification $\bar M\supset M$.
\end{proposition}

\section{Divisors on semi-tori.}

{\bf (a)}
Let $M$ be a complex semi-torus as before:
\begin{align*}
0 \to (\mathbf{C}^*)^p \to M&\, \overset{\eta}\to\, M_0 \to 0,\\
\bar\eta: \bar M& \to M_0.
\end{align*}
Let $D$ be an effective divisor on $M$ such that $D$ is compactified
to $\bar D$ in $\bar M$; that is, roughly speaking,
$D$ is algebraic along the fibers of $M \to M_0$.
If $M$ is a semi-Abelian variety, then this condition
is equivalent to the algebraicity of $D$.
We equip $L(\bar D)\to \bar M$ with a hermitian fiber metric.
Let $f:\mathbf{C} \to M$ be a holomorphic curve
such that $f(\C)\not\subset D$.
Let $\Omega$ be as in (\ref{3.4}).
Then there is a positive constant $C$ independent of $f$ such that
\begin{equation}
\label{4.1}
T_f(r; L(\bar D))=N(r; f^*D)+m_f(r; \bar D)+O(1)
\leqq C T_f(r; \Omega) +O(1).
\end{equation}

\begin{lemma}
\label{4.2}
Let $M$, $\bar M$, $M_0$ be as above.
Let $L\to \bar M$ be a line bundle on $\bar M$.
Then there exist a divisor $E$ with $\supp E \subset \partial M$
and a line bundle $L_0\to M_0$ such that
$L\cong L(E)\otimes \bar\eta^*L_0$ (in the sense of
bundle isomorphism or linear equivalence);
moreover,  such $L_0 \to M_0$ is uniquely determined (up to isomorphism).
\end{lemma}
\begin{proof}
Note that $\bar\eta:\bar M\to M_0$ is a topologically trivial
$\P_1(\C)^p$-bundle over $M_0$.
Hence by K\"unneth formula we have
\begin{equation}
\label{4.3}
H^2(\bar M,\Z)=H^2(\P_1(\C),\Z)^p\oplus H^2(M_0,\Z).
\end{equation}
Since the higher direct image sheaves
${\mathcal R}^q\eta_*{\mathcal O}, q \geqq 1$, vanish,
it follows that
$H^*(\bar M,\O)\cong H^*(M_0,\O)$.
We deduce that the Picard group $\Pic(\bar M)$ is 
generated
by
$\bar\eta^*\Pic(M_0)$
and the subgroup  of $\Pic(\bar M)$ generated by the irreducible
components of $\partial M=\bar M\setminus M$.
Thus for $\bar L \to \bar M$ there exists a
divisor $E$ with $\supp E \subset \partial M$ such that
$L \tensor L(-E) \in \bar\eta^*\Pic(M_0)$;
the assertion follows.
\end{proof}

We denote by
\begin{equation}
\label{4.4}
\StD=\{x \in M; x+D=D\}^0
\end{equation}
the identity component of those $x\in M$
which leaves $D$ invariant by translation.
The complex semi-subtorus $\StD$ (cf.\ [NW99])
is called the {\it stabilizer} of $D$.

\begin{lemma}
\label{4.5}
{\rm (i)}
Let $Z$ be a divisor on $\bar M$ such that $Z\cap M$ is effective.
Let $L_0 \in \Pic (M_0)$ such that
$L(Z)\tensor \bar\eta^*L_0^{-1} \cong L(E)$ with $\supp E \subset \partial M$.
Then $c_1(L_0)\geqq 0$.

{\rm (ii)}
Let $D$ be an effective divisor on $M$ with compactification
$\bar D$ as above.
Assume that $\StD=\{0\}$.
Then $\bar D$ is ample on $\bar M$.
\end{lemma}

\begin{proof}
(i)  Assume the contrary.
Recall that $M_0$ is a compact complex torus with universal covering
$\pi_0:\C^m\to M_0$. We may regard the Chern class $c_1(L_0)$ as 
bilinear form on the vector space $\C^m$.
Suppose that $c_1(L_0)$ is not semi-positive definite.
Let $v\in\C^m$ with $c_1(L_0)(v,v)<0$ and let $W$ denote the orthogonal
complement of $v$ (i.e.~$W=\{w \in \C^m:c_1(L_0)(v,w)=0\}$).
Let $\mu$ be a semi-positive skew-Hermitian form on $\C^m$ such that
$\mu(v,\cdot)\equiv 0$ and $\mu|_{W\times W}>0$.
Now consider the $(n-1,n-1)$-form $\omega$ on $\bar M$ given
by
\begin{equation}
\label{4.6}
\omega=\Omega^p\wedge\bar\eta^*\mu^{m-1}.
\end{equation}
By construction we have $\omega\wedge\bar\eta^*c_1(L_0) <0$.
Let $Z=Z'+Z''$ so that $Z'$ is effective and has no
component of $\partial M$, and $\supp Z'' \subset \partial M$.

By the Poincar\'e duality,
$$
\int_{\bar M}c_1(L(Z))\wedge\omega = \int_{Z} \omega.
$$
Since $\omega\wedge c_1(L(E))=0$, we  have
$$
\int_{\bar M}c_1(L(Z))\wedge\omega=
\int_{\bar M}\bar\eta^*c_1(L_0)\wedge\omega <0.
$$
On the other hand,
$$
\int_{Z}\omega = \int_{Z'}\omega + \int_{Z''} \omega.
$$
Note that
$\int_{Z'}\omega\geqq 0$,
because $Z'$ is effective and $\omega\geqq 0$,
and that $\int_{Z''} \omega=0$, because $\supp Z'' \subset \partial M$,
and $\Omega^p$ vanishes on $\partial M$ by construction.
Thus we deduced a contradiction.

(ii) When $p=0$, the assertion is well known ([We58]).
Assume $p>0$.
Let $\C^*$ act on $\bar M$ as the $k$-th factor of
$(\C^*)^p\subset M$. Since $\StD=\{0\}$, one infers
that there is an orbit whose closure intersects $\bar D$
transversally.
Hence,
\begin{equation}
\label{4.7}
c_1(L)= \left(n_1,\ldots,n_p;c_1(L_0)\right)
\end{equation}
in the form described in (\ref{4.3}) with $n_1,\ldots,n_p>0$.

Now let us consider $L_0$ as in the above (i).
By (i) we know that $c_1(L_0)\geqq 0$.
Assume that there is a vector $v\in\C^m\setminus \{O\}$ with
$c_1(L_0)(v,v)=0$.
Then we choose $\mu$ and $\omega$ as in (\ref{4.6}).
Because of the definition we have
\begin{equation}
\label{4.8}
\int_{\bar D} \omega = 0.
\end{equation}
By the flat connection of the bundle $\eta:M \to M_0$,
the vector $v$ is identified as a vector field on $M$.
Observe that $\bar\eta(\bar D)=M_0$.
The construction of $\omega$ and (\ref{4.8})
imply that $v\in T_x(D)$ for all $x\in D$.
It follows that the one-parameter subgroup corresponding
to $v$ must stabilize $D$;
this is a contradiction.
Thus $c_1(L_0)>0$ if $\StD=\{0\}$.

Since all $n_i>0$ in (\ref{4.7}) and $c_1(L_0)>0$,
it follows that $c_1(L(\bar D))$ is positive.
Thus $\bar D$ is ample on $\bar M$.
\end{proof}

\begin{corollary}
\label{4.9}
Let $f:\C \to M$ and $D$ be as above,
and let $\Omega$ be as in (\ref{3.4}).
Assume that $\StD=\{0\}$.
Then we have the following.
\begin{enumerate}
\item
There is a positive constant $C$ such that
$$
C^{-1}T_f(r; \Omega)+O(1)\leqq T_f(r; c_1(\bar D))\leqq
C T_f(r; \Omega)+O(1).
$$
\item
$S_f(r; \Omega)=S_f(r; c_1(\bar D))$.
\end{enumerate}
\end{corollary}

The proof is clear.

{\it Remark.}
$\bar D$ may be ample even if $\StD\ne\{0\}$.
For instance, this happens
for the diagonal divisor $D$ in $M=\C^*\times\C^*\hookrightarrow
\bar M=\P_1\times\P_1$.

{\bf (b) Boundary condition for $D$.}
We keep the previous notations.
Let
\begin{equation}
\label{4.10}
\partial M=\bigcup_{j=1}^p B_j
\end{equation}
be the Whitney stratification
of the boundary divisor of $M$ in $\bar M$;
that is, $B_j$ consists of all points $x \in \partial M$
such that the number of irreducible components of $\partial M$
passing $x$ is exactly $j$.
Set $B_0=M$.
A connected component of $B_j, 0\leqq j \leqq p$, is called
a {\it stratum} of the stratification $\bar M=\bigcup_{j=0}^pB_j$.
Observe that $\dim B_j=n-j$.

Note that the holomorphic action of $M$ on $M$ by translations
is equivariantly extended to an action on $\bar M$,
which preserves every stratum of $B_j, 0 \leqq j \leqq p$.

Let $D$ be an effective divisor of $M$
which can be extended to a divisor $\bar D$ on $\bar M$
by taking its topological closure of the support.
We consider the following boundary condition for $D$:

\begin{condition}
\label{4.11}
$\bar D$ does not contain any stratum of $B_p$.
\end{condition}
Note that the strata of $B_p$ are minimal.

\begin{lemma}
\label{4.12}
If condition \ref{4.11} is fulfilled,
then
$$
\dim \bar D \cap B_j < \dim B_j= n - j,
\quad 0 \leqq \forall j \leqq p.
$$
\end{lemma}
\begin{proof}
Assume the contrary. Then there exists a stratum $S\subset B_j$
such that $S\subset \bar D$. Clearly the closure $\bar S$ of $S$ 
is likewise contained in $\bar D$. But the closure of any
stratum contains a minimal stratum, i.e., contains a stratum
of $B_p$. However, this is in contradiction
to condition \ref{4.11}.
\end{proof}
\begin{example}
\label{4.13}{\rm\hskip5pt
Take a classical case where $M$ is the complement of $n+1$
hyperplanes $H_j$ of $\P^n(\C)$ in general position.
Then $M\cong(\C^*)^n$.
Let $D=H_{n+2}$ be an $(n+2)^{\mathrm{th}}$ hyperplane of $\P^n(\C)$.
Then condition \ref{4.11} is equivalent to that
all $H_j, 1 \leqq j \leqq n+2 $, are in general position.}
\end{example}

Next we interpret the boundary condition \ref{4.11}
in terms of local defining equations of $\bar D$.
Take $\sigma\in H^0(\bar M,L(\bar D))$ such that $(\sigma)=\bar D$.
Suppose that $p >0$.
Let $x_0 \in \partial M \cap \bar D$ be an arbitrary point.
Let $E$ and $L_0$ be as in Lemma \ref{4.2} for $L=L(\bar D)$.
We take an open neighborhood $U$ of $\bar\eta(x_0)$ such that
the restrictions $\bar M|U$ and $L_0|U$ to $U$ are trivialized.
Write
$$
x_0=(u_0,x_0'')\in (\P^1(\C))^p\times U\cong\bar M|U.
$$
We take an open neighborhood $V$ of $u_0$ such that
$V \cong \C^p \subset (\P^1(\C))^p$ with coordinates
$(u_1,\ldots,u_p)$.
Then $L(\bar D)|(V\times U)$ is trivial, and hence
$\sigma|(V\times U)$ is given by a polynomial function
\begin{equation}
\label{4.14}
\sigma(u, x'')= \sum_{\mathrm{finite}}
a_{l_1 \cdots l_p}(x'')
u_1^{l_1} \cdots u_p^{l_p},
\quad (u,x'') \in V\times U,
\end{equation}
with coefficients $a_{l_1 \cdots l_p}(x'')$ holomorphic in $U$.
Since $\bar D$ has no component of $\partial M$,
$\sigma(u, x'')$ is not divisible by any $u_j$.
Set $u_0=(u_{01}, \ldots, u_{0p})$.
Then, after a change of indices of $u_i$ one may assume that
$u_{01}= \cdots =u_{0q}=0, u_{0i} \not= 0, 1\leqq q < i \leqq p$.
Expand $\sigma(u, x'')$ and set $\sigma_1$ and $\sigma_2$ as follows:
\begin{align}
\label{4.15}
\sigma(u, x'') &= \sum_{l_1+\cdots + l_q \geqq 1}
a_{l_1 \cdots l_p}(x'')
u_1^{l_1} \cdots u_p^{l_p}+
\sum_{l_1=\cdots = l_q=0}
a_{0 \cdots 0 l_{q+1} \cdots l_p}(x'')
u_{q+1}^{l_1} \cdots u_p^{l_p}, \\
\sigma_1 &= \sum_{l_1+\cdots + l_q \geqq 1}
a_{l_1 \cdots l_p}(x'')
u_1^{l_1} \cdots u_p^{l_p}, \nonumber\\
\sigma_2 &=\sum_{l_1=\cdots = l_q=0}
a_{0 \cdots 0 l_{q+1} \cdots l_p}(x'')
u_{q+1}^{l_1} \cdots u_p^{l_p} . \nonumber
\end{align}
We have

\begin{lemma}
\label{4.16}
Let the notation be as above.
Then condition \ref{4.11} is equivalent to that
for every $x_0 \in \partial M$, $\sigma_2\not\equiv 0$.
\end{lemma}

{\bf (c) Regularity of stabilizers.}
Let $M$ be a complex semi-torus with fixed presentation
as in (\ref{3.1}):
\begin{equation}
\label{4.17}
0 \to G=(\C^*)^p \to M \to M_0 \to 0.
\end{equation}

{\it Definition.}
A closed complex Lie subgroup $H$ of $M$ is called {\it regular}
if there is a subset $I\subset\{1,\ldots,p\}$ such that
$$
G \cap H=\{(z_1,\ldots,z_p)\in G ; z_i=1, \forall i\in I\}.
$$

Regular subgroups are those compatible with the compactification
induced by (\ref{4.17}).
The presentation (\ref{4.17}) induces in a canonical way such
presentations for $H$ and $M/H$.

\begin{lemma}
\label{4.18}
Let $H$ be a regular Lie subgroup of $M$.
Then the quotient mapping $M\to M/H$ is extended holomorphically
in  a natural way to the compactification
$$
\bar M \stackrel{\bar H}{\longrightarrow}\overline{(M/H)},
$$
which is a holomorphic fiber bundle of compact complex manifolds
with fiber $\bar H$.
\end{lemma}

We will prove the following proposition.

\begin{proposition}
\label{4.19}
Let $M$ be a semi-torus with presentation (\ref{4.17}) and
let $D$ be an effective divisor fulfilling
the boundary condition \ref{4.11}.
Then there exists a finite unramified covering $\mu':M_0'\to M_0$
such that $\mathrm{St}({\mu}^* D)$
is regular in $M'$,
where $\mu:M'\to M$ is the finite covering of $M$ induced by $\mu'$;
i.e., $M'=M\times_{M_0}M_0'$.
\end{proposition}

{\it Remark.}
Note that $\mu$ extends holomorphically to the unramified covering
of the compactification $\bar M$,
$\bar\mu :\bar M' \to \bar M$.
\medskip

\begin{proof}
First, if $D$ is invariant under one of the $p$ direct factors of
$G=(\C^*)^p$ in (\ref{4.17}),
we take the corresponding quotient.
Thus we may assume that $\StD\cap G$ does not contain anyone
of the $p$ coordinate factors of $G$.

Assume that $\dim \StD \cap G >0$.
Let $I$ be a subgroup of $\StD\cap G$ isomorphic to $\C^*$.
Then there are integers $n_1,\ldots,n_p$ such that
\[
I=\{(t^{n_1},\ldots,t^{n_p}):t\in\C^*\}.
\]
By re-arranging indices and coordinate changes
of type, $z_i\mapsto \frac{1}{z_i}$,
we may assume that there is a natural number $q$ such that
$n_i>0$ for $i\leqq q$ and $n_i=0$ for $i>q$.
Let $G=G_1\times G_2$ with 
\begin{align*}
G_1 &=\{(u_1,\ldots,u_{q-1},\underbrace{1, \ldots,1}_{p-q});
u_i\in\C^*\} \subset G,\\
G_2 &=\{(\underbrace{1, \ldots,1}_{q},
u_{q+1},\ldots,z_p); u_i \in \C^* \} \subset G.
\end{align*}
Then $I\subset G_1$.
Consider $\lambda :M \to M/G_1$.
If $\lambda(D) \ne M/G_1$, then $D$ would be
$G_1$-invariant and in particular would be invariant under
the coordinate factor groups contained in $G_1$.
Since this was ruled out, we have $\lambda(D)=M/G_1$.
Now observe that for every
$u=(u_1,\ldots,u_p)\in\C^p \subset (\P^1(\C))^p$ we have
$$
\lim_{t\to 0} \; (t^{n_1},\ldots,t^{n_p})\cdot u
=(0,\ldots,0,u_{q+1},u_{q+2},\ldots, u_p).
$$
Hence it follows from $I\subset\StD$ and $\lambda(D)=M/G_1$ that
$$
\{0\}^q \times (\P^1(\C))^{p-q} \subset \bar D.
$$
This violates the boundary condition \ref{4.11} because of Lemma \ref{4.12}.
Thus $G\cap\StD$ is zero-dimensional, and hence finite.
As a consequence, $\StD$ is compact.
After a finite covering, $\StD$ maps injectively in $M_0$ and
therefore is regular.
\end{proof}

\section{Proof of the Main Theorem.}

We first prove the following key lemma:
\begin{lemma}
\label{5.1}
Assume the same conditions as in the Main Theorem.
Then,
$$
m_f(r; \bar D)=S_f(r; c_1(\bar D)).
$$
\end{lemma}

Besides the conditions stated above,
we may also assume by Proposition \ref{4.19}
and Lemma \ref{4.5}, (ii) that $\StD=\{0\}$, $\bar D$ is ample on $\bar M$,
and hence $M$ is a semi-Abelian variety $A$:
$$
0 \to (\C^*)^p \to A \to A_0 \to 0.
$$
We keep these throughout in this section.

Here we need the notion of logarithmic jet spaces due to [No86].
Since $\partial A$ has only normal crossings,
we have the logarithmic $k$-th jet bundle
$J_k(\bar A; \log \partial A)$ over $\bar A$ along
$\partial A$, and a morphism
$$
\psi_k:J_k(\bar A; \log \partial A) \to J_k(\bar A)
$$
such that the sheaf of germs of holomorphic sections of
$J_k(\bar A; \log \partial A)$ is isomorphic to
that of logarithmic $k$-jet fields (see [No86], Proposition (1.15);
there, a ``subbundle'' $J_k(\bar A; \log \partial A)$ of
$J_k(\bar A)$ should be understood in this way).
Because of the flat structure of the logarithmic tangent
bundle $\mathbf{T}(\bar A; \log \partial A)$,
$$
J_k(\bar A; \log \partial A) \cong \bar A
\times \mathbf{C}^{nk}.
$$
Let
\begin{align}
\label{5.2}
\pi_1 &: J_k(\bar A; \log \partial A)
\cong \bar A \times \mathbf{C}^{nk} \to \bar A, \\
\pi_2 &: J_k(\bar A; \log \partial A)
\cong \bar A \times \mathbf{C}^{nk} \to \mathbf{C}^{nk}\nonumber
\end{align}
be the first and the second projections.
For a $k$-jet $y \in J_k(\bar A; \log \partial A)$ we call
$\pi_2(y)$ the {\it jet part} of $y$.

Let $x \in \bar D$ and let $\sigma=0$ be a local defining
equation of $\bar D$ about $x$.
For a germ $g:(\mathbf{C}, 0) \to (A,x)$ of
a holomorphic mapping we denote its $k$-jet by $j_k(g)$ and
write
$$
d^j\sigma (g)= \left. \frac{d^j}{d\zeta^j}\right|_{\zeta=0}
\sigma(g(\zeta)).
$$
We set
\begin{align*}
J_k(\bar D)_x &= \left\{j_k(g)\in J_k(\bar A)_x; d^j\sigma(g)=0,
1 \leqq j \leqq k \right\},\\
J_k(\bar D) &= \bigcup_{x \in D} J_k(\bar D)_x ,\\
J_k(\bar D;\log\partial A)&= \psi_k^{-1}J_k(\bar D).
\end{align*}
Then $J_k(\bar D; \log\partial A)$ is a subspace of
$J_k(\bar A;\log\partial A)$, which is depending in general on the
embedding $\bar D \hookrightarrow \bar A$
(cf.\ [No86]).
Note that $\pi_2(J_k(\bar D; \log \partial A))$
is an algebraic subset of $\mathbf{C}^{nk}$.

Let $J_k(f):\mathbf{C} \to J_k(\bar A; \log\partial A)=
\bar A\times\mathbf{C}^{nk}$
be the $k$-th jet lifting of $f$.
Then by [No98] the Zariski closure of
$J_k(f)(\mathbf{C})$ in $J_k(\bar A; \log\partial A)$ is of the form,
$\bar A\times W_k$, with an affine irreducible subvariety
$W_k \subset \mathbf{C}^{nk}$.
Let $\pi:\mathbf{C}^n \to A$ be the universal covering
and let
$$
\tilde f:z \in \mathbf{C}\to
\left(\tilde f_1(z), \ldots, \tilde f_n(z)\right)
\in \mathbf{C}^n
$$
be the lifting of $f$.
Assume that $f$ is of finite order.
Then $\tilde f(z)$ is a vector valued polynomial by Lemma \ref{3.9}.
Note that every non-constant polynomial map from $\C$ to $\C^n$ is
proper, and hence the image is an algebraic subset.
It follows that
$$
W_k=\overline{\left\{\left(\tilde f'(z), \ldots, \tilde f^{(k)}(z)\right)
z \in \mathbf{C}\right\}}
=\left\{\left(\tilde f'(z), \ldots, \tilde f^{(k)}(z)\right),
z \in \mathbf{C}\right\},
$$
and hence $\dim W_k \leqq 1$.
Thus we deduced the following lemma.

\begin{lemma}
\label{5.3}
Let the notation be as above.
If $f:\C\to A$ is of finite order,
then $\dim W_k \leqq 1$ and 
for every point $w_k \in W_k$ there is a point $a \in \mathbf{C}$
with $\pi_2 \circ J_k(f)(a)=w_k$.
\end{lemma}

\begin{lemma}
\label{5.4}
Let the notation be as above.
\begin{enumerate}
\item
Suppose that $f$ is of finite order $\rho_f$.
Then there is a number $k_0=k_0(\rho_f,D)$ such that
$$
\pi_2(J_k(\bar D;\log\partial A)) \cap W_k = \emptyset, \quad k \geqq k_0.
$$
\item
Suppose that $f$ is of infinite order.
Then there is a number $k_0=k_0(f,D)$ such that
$$
\pi_2(J_k(\bar D;\log\partial A)) \cap W_k \not= W_k, \quad k \geqq k_0.
$$
\end{enumerate}
\end{lemma}

{\it Proof.}
(i)  By making use of (\ref{5.2}) we have the projection
$p_{k,l}:\mathbf{C}^{nk} \to \mathbf{C}^{nl}$ for $k \geqq l$
induced from the canonical projection
$J_k(\bar A; \log \partial A) \to J_l(\bar A; \log \partial A)$.
For a subset or a point $E_k$ of $\mathbf{C}^{nk}$ and $l \leqq k$
we write $E_{k,l}=p_{k,l}(E_k)$.

We see first by Lemma \ref{3.9} that $\rho_f \in \mathbf{Z}$,
and $\tilde f (z)$ is a vector valued
polynomial of order $\leqq \rho_f$.
Thus, $W_k$ is of form
$$
W_k=(W_{k,\rho_f},\; \underbrace{O, \ldots, O}_{k-\rho_f}\, ).
$$
Set $W'_k=W_k \cap \pi_2(J_k(\bar D; \log \partial A))$.
Then we have
$$
W'_k=(W'_{k,\rho_f},\; \underbrace{O, \ldots, O}_{k-\rho_f}\, ).
$$
Assume that the present assertion fails.
Then, by the Noetherian property of algebraic subsets,
there is a point $\xi_{\rho_f} \in \bigcap_{k=\rho_f}^\infty W'_{k,\rho_f}$
such that, setting
$\xi_k=(\xi_{\rho_f},\; \underbrace{O,\ldots, O}_{k-\rho_f}\,)
\in \mathbf{C}^{nk}$,
we have
$$
\xi_k \in \pi_2(J_k(\bar D; \log\partial A)), \quad \forall k \geqq \rho_f.
$$
We identify $\xi_k$ with a logarithmic $k$-jet field on
$\bar A$ along $\partial A$ (see [No86]).
Set $S_k=\pi_1\big(J_k(\bar D; \log\partial A) \cap \pi_2^{-1}
(\xi_k)\big)$.
Then,
$$
\bar D \supset S_{\rho_f} \supset S_{\rho_f+1}\supset \cdots,
$$
which stabilize to $S_0=\bigcap_{k=\rho_f}^\infty S_k \not=\emptyset$.
Let $x_0 \in S_0.$
If $x_0 \in A$, it follows from \ref{5.3} that
there are points $a \in \mathbf{C}$ and $y_0 \in A$ such that
\begin{align*}
f(a)+x_0+y_0 &\in D, \\
\left. \frac{d^k}{dz^k}\right|_{z=a} \sigma(f(z)) &=0, \quad
\forall k \geqq 1,
\end{align*}
where $\sigma$ is a local defining function of $D$
about $f(a)+x_0+y_0$.
Therefore
$$
f(\mathbf{C})+x_0+y_0 \subset D,
$$
and hence this contradicts the Zariski denseness
of $f(\mathbf{C})$ in $A$.
This finishes the proof in the case of $x_0 \in A$.

Suppose now that $x_0 \in \bar A\setminus A$.
Let $\partial A=\bigcup B_j$ be the Whitney stratification
as in (\ref{4.10}), and let $x_0 \in B_q$.
Let $B$ be the stratum of $B_q$ containing $x_0$.
Then $B$ itself is a semi-Abelian variety such that
$$
0 \to (\C^*)^{p-q} \to B \to A_0 \to 0.
$$
Let $\sigma(u, x'')=\sigma_1(u, x'')+\sigma_2(u, x'')$
be as in (\ref{4.15}) and define $\bar D$ in a neighborhood
$W$ of $x_0$ such that $W$ is of type $V\times U$ as in (\ref{4.14}).
It follows from Lemma \ref{4.16} that $\sigma_2\not\equiv 0$.
Note that $\bar D \cap W\cap B$ is defined by $\sigma_2 = 0$ in $B$.
There is a point $a \in \mathbf{C}$ such that
$\pi_2\circ J_{\rho_f}(f)(a)=\xi_{\rho_f}$.
Dividing the coordinates into three blocks, we set
$$
x_0=(\underbrace{0, \ldots, 0}_{q},\; x'_{0},\; x''_0).
$$
We may regard $w_0=(x'_{0},\;  x''_0) \in B$.
Taking a shift $f(z)+y_0$ with $y_0 \in A$ so that
$f(a)+y_0 \in W$, we set in a neighborhood of $a \in \C$
\begin{align}
\label{5.5}
f(z)+y_0 &=(u_1(z), \ldots, u_q(z),\: u_{q+1}(z), \ldots,
u_p(z),\: x''(z)) \in W, \\
g(z) &=(u_{q+1}(z), \ldots, u_p(z),\: x''(z)) \in W\cap B. \nonumber
\end{align}
Here we may choose $y_0$ so that $g(a)=w_0$.

We set $\xi_k=\pi_2 \circ J_k(f)(a)$ for all $k \geqq 1$.
Using the same coordinate blocks as (\ref{5.5}), we set
\begin{align*}
\xi_k &= (\xi'_{k(1)}, \xi'_{k(2)}, \xi''_k),\qquad k \geqq \rho_f, \\
\xi_{k(2)} &= \hbox{ the jet part of }J_k(g)(a) =(\xi'_{k(2)}, \xi''_k).
\end{align*}
Since the logarithmic term (e.g.,
$z_j\frac{\partial}{\partial z_j},1 \leqq j \leqq q$,
in the case of 1-jets) of a logarithmic jet field
vanishes on the corresponding divisor locus
(e.g., $\bigcup_{j=1}^q\{z_j=0\}$)
(see [No86], \S1 and (1.14) for more details),
we have $\xi_{k}(\sigma_1)(x_0)=0$ by (\ref{4.15}),
and hence $\xi_{k(2)}(\sigma_2)(x_0)=0, \forall k\geqq 1$; i.e.,
\begin{equation}
\label{5.6}
\left. \frac{d^k}{dz^k}\right|_{z=a} \sigma_2(g(z)) =0, \quad
\forall k \geqq 0.
\end{equation}
Let $(\mathbf{C}^*)^q$ be the first $q$-factor of the
subgroup $(\mathbf{C}^*)^p \subset A$, and let
$\lambda: A \to A/(\mathbf{C}^*)^q \cong B$ be the quotient map.
By (\ref{5.5}) and (\ref{5.6}) the composed map
$g(z)=\lambda \circ (f(z)+y_0)$ has
an image contained in $\bar D \cap B$; that is,
it has no Zariski dense image in $A/(\mathbf{C}^*)^q$,
and hence so is $f$;
this is a contradiction.

The order of the tangency of $f$ and the above used $g$ with
$\bar D$ is bounded
as $\tilde f$ runs over all vector valued
polynomials of order at most $\rho_f$ such that $f(\mathbf{C})$
is Zariski dense in $A$.
Hence there is such a number $k_0$ depending only on $\rho_f$
and $D$.

(ii) Assume contrarily that
$\pi_2(J_k(\bar D; \log \partial A)) \cap W_k = W_k$
for all $k\geqq 1$.
Since $\pi_2 \circ J_k(f)(0) \in W_k,  \forall k \geqq 1$,
we apply the same argument as in (i) with setting
$\xi_k=\pi_2 \circ J_k(f)(0)$.
Then we deduce a contradiction that $f$ has no Zariski dense image.
{\it Q.E.D.}
\medskip

{\it Proof of Lemma \ref{5.1}.}
For a multiple $l \bar D$ of $\bar D$ we have
$$
m_f(r; l\bar D)=lm_f(r;\bar D).
$$
Thus we may assume that $\bar D$ is very ample on $\bar A$.
Let $\{\tau_j\}_{j=1}^N$ be a base of $H^0(\bar A, L(\bar D))$
such that $\supp (\tau_j)\not\supset f(\C)$ for all $1\leqq j\leqq N$.
Since $\bar D$ is very ample, the sections $\tau_j$, $1\leqq j\leqq N$,
have no common zero.
Set
$$
U_j =\{\tau_j\not=0\}, \quad 1 \leqq j \leqq N.
$$
Then $\{U_j\}$ is an affine open covering of $\bar A$.
Let $\sigma \in H^0(\bar A, L(\bar D))$ be a section such that
$(\sigma)=\bar D$.
We define a regular function $\sigma_j$ on every $U_j$ by
$$
\sigma_j(x)=\frac{\sigma(x)}{\tau_j(x)}.
$$
Note that $\sigma_j$ is a defining function of $\bar D \cap U_j$.
Let us now fix a hermitian metric $||\cdot||$ on $L(\bar D)$.
Then there are positive smooth functions $h_j$ on $U_j$
such that
\begin{equation}
\nonumber
\frac{1}{\|\sigma (x)\|}=\frac{h_j(x)}{|\sigma_j(x)|}, \quad x \in U_j.
\end{equation}

Assume that $f$ is of finite order.
By Lemma \ref{5.4} there are regular functions
$b_{ji}, 0 \leqq i \leqq k_0$, on $U_j \times W_{k_0}$
such that
\begin{equation}
\label{5.7}
b_{j0}\sigma_j+
b_{j1}d\sigma_j +\cdots+
b_{j k_0}d^{k_0}\sigma_{j}=1.
\end{equation}
Here every $b_{ji}$ is expressed as
$$
b_{ji}= \sum_{\mathrm{finite}}b_{jil \beta_l}(x) w_l^{\beta_l},
$$
where $b_{{i\alpha} j l\beta_l}(x)$ are regular
functions on $U_j$ and $w_l$ are restrictions
of coordinate functions of $\mathbf{C}^{nk_0}$
to $W_{k_0}$.
Thus we infer that in every $U_j$
\begin{equation}
\label{5.8}
\frac{1}{\| \sigma \|}=
\frac{h_j}{|\sigma_j |}=
\left|
h_j b_{j 0} +
h_j b_{j 1}
\frac{d\sigma_{j}}{\sigma_{j}}
+\cdots+h_j
b_{j k_0}
\frac{d^{k_0}\sigma_{j}}{\sigma_{j}}
\right|.
\end{equation}

Take relatively compact open subsets $U'_j \Subset U_j$
(in the sense of differential topology) so that
$\bigcup U'_j=\bar A$.
For every $j$ there is a positive constant $C_j$
such that for $x \in U'_j$
\begin{equation}
\nonumber
h_j|b_{ji}|
\leqq \sum_{\mathrm{finite}}h_j
|b_{j i l \beta_l}(x)|\cdot
 |w_l|^{\beta_l}
\leqq C_{j} \sum_{\mathrm{finite}}|w_l|^{\beta_l}.
\end{equation}
Thus, after making $C_j$ larger if necessary,
there is a number $d_{j}>0$ such that for
$f(z) \in U'_j$
\begin{equation}
\nonumber
h_j(f (z))|b_{ji}(J_{k_0}(f)(z))|
\leqq C_j\left(1+
\sum_{1\leqq l \leqq n,1\leqq k \leqq k_0}
\left| \tilde f^{(k)}_l (z) \right|
\right)^{d_{j}}.
\end{equation}
We deduce that
\begin{align}
\nonumber
\frac{1}{\|\sigma(f(z))\|} \leqq& \sum_{j=1}^N
 C_j\left( 1+
\sum_{1\leqq l \leqq n, 1\leqq k \leqq k_0}
\left| \tilde f^{(k)}_l (z) \right|
\right)^{d_{j}}\\
&\times
\left(1+
\left|\frac{d\sigma_j}{\sigma_j}(J_1(f)(z))\right|+\cdots+
\left|\frac{d^{k_0}\sigma_j}{\sigma_j}(J_{k_0}(f)(z))\right|
\right).\nonumber
\end{align}
Hence one gets
\begin{align}
\label{5.9}
m_f(r; \bar D) &= \frac{1}{2\pi} \int_0^{2\pi} \log^+
\frac{1}{\|\sigma(f(re^{i\theta}))\|}d\theta \\
&\leqq O
\bigg(
\sum_{1\leqq l \leqq n, 1 \leqq k \leqq k_0}
\frac{1}{2\pi}\int_0^{2\pi}\log^+
\left|
\tilde f_l^{(k)} (re^{i\theta})
\right|
d\theta
\nonumber\\
&\quad
+ \sum_{1\leqq j \leqq N, 1\leqq k \leqq k_0}
\frac{1}{2\pi}\int_0^{2\pi}\log^+
\left| 
\frac{d^k\sigma_j}{\sigma_j}
(J_k(f)(re^{i\theta}))
\right| d\theta
\bigg)+O(1) . \nonumber
\end{align}
Recall that the rational functions $\sigma_j$ are equal
to quotients of two holomorphic sections $\sigma$ and $\tau_j$
of $L(\bar D)$.
By Lemma \ref{2.5}, (ii) we see that
$$
\frac{1}{2\pi}\int_0^{2\pi}\log^+\left| 
\frac{d^k\sigma_j}{\sigma_j}
(J_k(f)(re^{i\theta}))
\right| d\theta =m\left(r, \frac{(\sigma_j\circ f)^{(k)}}
{\sigma_j \circ  f}\right)=O(\log r).
$$
This combined with (\ref{5.9}) and Lemma \ref{3.8}
implies that $m_f(r; \bar D)=O(\log r)$;
this completes the proof in the case of finite order.

Assume that $f$ is of infinite order.
It follows from Lemma \ref{5.4}, (ii) that there exists a
polynomial function $R(w)$ in $w \in W_{k_0}$ such that
$$
\pi_2(J_k(\bar D;\log\partial A)) \cap W_k \subset
\{w \in W_{k_0}; R(w)=0\} \not=W_{k_0}.
$$
We regard $R$ as a regular function on every
$U_j\times W_{k_0}$.
Then we have the following equation
on every $U_j\times W_{k_0}$ with coefficients
similar to those of (\ref{5.7}):
\begin{equation}
\label{5.10}
b_{j 0}\sigma_{j}+
b_{j 1}d\sigma_{j}+\cdots+
b_{j k_0}d^{k_0}\sigma_{j}=R.
\end{equation}
Then, after the same arguments as in the case of finite order,
we have that for $f(z) \in U'_j$
\begin{align}
\label{5.11}
\frac{1}{\|\sigma(f(z))\|} =&
\frac{1}{\left|R\left(\tilde f^\prime(z),
\ldots, \tilde f^{(k_0)}(z)\right)\right|}
 \\
&\times
\left|
h_j b_{j 0} + h_j b_{j 1}
\frac{d\sigma_{j}}{\sigma_{j}}
+\cdots+h_j
b_{j k_0}
\frac{d^{k_0}\sigma_{j}}{\sigma_{j}}
\right| \nonumber \\
\leqq&
\frac{1}{\left|R\left(\tilde f^\prime(z), \ldots,
\tilde f^{(k_0)}(z)\right)\right|}
\sum_{j'=1}^N C_{j'}\left(1+
\sum_{1\leqq l \leqq n,1\leqq k \leqq k_0}
\left| \tilde f^{(k)}_l (z) \right|
\right)^{d_{j'}}
 \nonumber\\
&\times
\left(1+
\left|\frac{d\sigma_{j'}}{\sigma_{j'}}(J_1(f)(z))\right|+\cdots+
\left|\frac{d^{k_0}\sigma_{j'}}{\sigma_{j'}}(J_{k_0}(f) (z))\right|
\right) .
\nonumber
\end{align}
It follows that
\begin{align}
\nonumber
m_f(r; \bar D) \leqq& \frac{1}{2\pi}\int_0^{2\pi} \log^+
\frac{1}{\|\sigma(f(re^{i\theta}))\|}d\theta +O(1)\\
\leqq&
m\left(r,
\frac{1}{R\left(\tilde f^\prime, \ldots, \tilde f^{(k_0)}\right)}
\right)
+O\left(
\sum_{1\leqq l \leqq n,1 \leqq k \leqq k_0}
m\left(r, \tilde f_l^{(k)}\right) \right. \nonumber\\
&\left. + \sum_{1 \leqq j \leqq N, 1 \leqq k \leqq k_0}
m\left(r, \frac{d^k\sigma_j}{\sigma_j}\circ J_k(f)\right)
\right)+O(1)
\nonumber\\
\leqq&
T\left(r,
R\left(\tilde f^\prime, \ldots, \tilde f^{(k_0)}\right)\right)
\nonumber\\
&+O
\left(
\sum_{l, k, j}
m(r, \tilde f_j^{(k)})
 + m\left(r, \frac{(\sigma_j\circ f)^{(k)}}{\sigma_j\circ f}
\right)\right) +O(1).
\nonumber
\end{align}
This combined with Lemmas \ref{3.8} and \ref{2.5} implies that
$m_f(r; \bar D)=S_f(r; c_1(\bar D))$.
This finishes the proof.
{\it Q.E.D.}
\medskip

{\it Proof of the Main Theorem.}
We keep the notation used above.
Thanks to Lemma \ref{5.1} the only things we still
have to show are the statements on the truncation, i.e.,
the bounds on $N(r;f^*D)- N_{k_0}(r;f^*D)$.
Observe that
$\ord_z f^*D>k$ if and only if $J_{k}(f)(z)\in J_k(\bar D; \log \partial A))$.
Therefore, if $f$ is of finite order,
Lemma \ref{5.4}, (i) implies that $N(r;f^*D)=N_{k_0}(r;f^*D)$.

In the case where $f$ is of infinite order
we infer from (\ref{5.11}) that
$$
\ord_z f^*D - \min\{\ord_z f^*D, k_0\} \leqq
\ord_z \left(R\left(\tilde f', \ldots, \tilde f^{k_0}\right)\right)_0.
$$
Thus we have after integration that
$$
N(r; f^*D) - N_{k_0}(r; f^*D) \leqq
N\left(r; \left(R\left(\tilde f^\prime, \ldots, \tilde
f^{(k_0)}\right)\right)_0\right).
$$
It follows from (\ref{2.2}), (\ref{2.3}) and Lemma \ref{3.8} that
\begin{align*}
N\left(r; \left(R\left(\tilde f^\prime, \ldots, \tilde
f^{(k_0)}\right)\right)_0\right)
&\leqq 
T\left(r, R\left(\tilde f^\prime, \ldots, \tilde
f^{(k_0)}\right)\right)+O(1)\\
&\leqq O\left(
\sum_{1\leqq l \leqq n, 1\leqq k \leqq k_0} T\left(r, \tilde f_l^{(k)}\right)
\right)\\
&= S_f(r; \Omega).
\end{align*}
Furthermore, $S_f(r; \Omega)=S_f(r,c_1(D))$ by Corollary \ref{4.9}, (ii),
because $\bar D$ is ample.
Hence,
$$
N(r; f^*D) \leqq N_{k_0}(r; f^*D) + S_f(r; c_1(\bar D)).
$$
The proof is completed.
{\it Q.E.D.}
\medskip

By (\ref{3.5}), (\ref{3.9}) and Lemma \ref{5.1} we have

\begin{corollary}
\label{5.12}
Let $M$ be a complex torus and let $f: \C \to M$ be an
arbitrary holomorphic curve.
Let $D$ be an effective divisor on $M$ such that
$D \not\supset f(\C)$.
Then we have the following.
\begin{enumerate}
\item
Suppose that $f$ is of finite order $\rho_f$.
Then there is a positive integer $k_0=k_0(\rho_f, D)$ such that
$$
T_f(r;c_1(D))=N_{k_0}(r; f^*D) + O(\log r).
$$
\item
Suppose that $f$ is of infinite order.
Then there is a positive integer $k_0=k_0(f, D)$ such that
$$
T_f(r;c_1(D))=N_{k_0}(r; f^*D) + S_f(r; c_1(D)).
$$
\end{enumerate}
Specially, $\delta(f; D)=\delta_{k_0}(f; D)=0$ in both cases.
\end{corollary}

{\it Proof.}
Since the Zariski closure of $f(\C)$ is a translation
of a complex subtorus of $M$
(cf., e.g., [\NO], Chap.\ VI, [Ko98], Chap.\ 3, \S9, [NW99]),
we  may assume that $f(\C)$ is Zariski dense.
Hence this statement is a special case of the Main Theorem.
{\it Q.E.D.}
\medskip

\begin{proposition}
\label{5.13}
Let $M$ be a complex semi-torus $M$ and let $D$ 
be an effective divisor on $M$
such that its topological closure $\bar D$ is a divisor in
$\bar M$.
Assume that $D$ violates the boundary condition \ref{4.11}.
Then there exists an entire holomorphic curve
$f:\C\to M$ of an arbitrarily given integral order
$\rho\geqq 2$ in general, and $\rho\geqq 1$
in the case of $M_0=\{0\}$ such that
$f(\C)$ is Zariski dense in $M$ and
$\delta(f;\bar D)>0$.
\end{proposition}

\begin{proof}
Let $\hat M=(\P^1(\C))^p\times\C^m \to \bar M$
(resp.\ $\C^m \to M_0$)
be the universal covering of $\bar M$ (resp.\ $M_0$),
and $\hat D \subset \hat M$ the preimage of $\bar D$.
We may assume that
\[
\{(\infty)\}^p\times \C^m \subset \hat D.
\]
Let $c_1,\ldots,c_p$ be $\Q$-linear independent real numbers
with
\begin{equation}
\label{5.14}
0 < c_1 < c_2 < \cdots < c_p.
\end{equation}
Let $\rho\geqq 2$ or $\rho\geqq 1$ be an arbitrary integer
as assumed in the proposition, and set
\begin{equation}
\label{5.15}
\hat f:z\mapsto \left(\left[1:e^{c_1z^\rho}\right],
\left[1:e^{c_2z^\rho}\right],\ldots,
\left[1:e^{c_pz^\rho}\right];L(z)\right),
\end{equation}
where $L:\C\to\C^m$ is a linear map such that the image $L(\C)$
in $M_0$ is Zariski dense.
Moreover, by a generic choice of $c_j$ and $L$ we have that
$f(\C)$ is Zariski dense in $M$.
Let \[ U_i\Subset V_i\Subset M_0
\]
be a finite collection of relatively compact holomorphically convex
open subsets of $M_0$ such that there are sections
$\mu_i:V_i \overset{\scriptstyle\sim}{\to} \hat V_i\subset\C^m$ and
such that the $U_i$ cover $M_0$. Set $\hat U_i=\mu_i(U_i)$.

For every $i$ the restricted divisor
$\hat D|((\P^1(\C))^p\times \hat V_i)$
is defined by a homogeneous polynomial $P_{i0}$ of multidegree
$(d_1,\ldots,d_p)$, where the coefficients are holomorphic
functions on $V_i$.
Let $P_i$ denote the associated inhomogeneous polynomial.
Then $P_i$ is a polynomial of multidegree $(d_1,\ldots,d_p)$.
Due to $\{\infty\}^p\times\C^m\subset\hat D$, $P_i$ does not carry
the highest degree monomial, $u_1^{d_1}\cdots u_p^{d_p}$.

Recall $\bar M=\hat M/\Lambda_0$ where $\Lambda_0$
is a lattice in $\C^m$
and acts on $\hat M$ via 
\[
\lambda:(u_1,\ldots,u_p; x'') \mapsto
\lambda\cdot (u; x'')=
(\beta_1(\lambda)u_1,\ldots,\beta_p(\lambda)u_p; x''+\lambda),
\]
where $\beta:\Lambda_0 \to (S^1)^p$ is a group homomorphism
into the product of $S^1=\{|z|=1; z \in \C^* \}$.

Together with (\ref{5.15}) and (\ref{5.14}) this implies that
there is a constant $C>0$ such that
\begin{equation}
\label{5.16}
|P_i(\lambda\cdot\hat f(z)) |\leqq 
C |e^{(\sum_j d_jc_j)z^\rho-c_1z^\rho}|
\end{equation}
for all $\lambda\in\Lambda_0$ and $z\in\C$ with $\re z^\rho >0$
and $\lambda\cdot\hat f(z)\in(\P^1(\C))^p\times {\hat U_i}$.
Note that for every $z\in\C$ there exists an element $\lambda\in\Lambda_0$
and an index $i$ such that
$\lambda\cdot\hat f(z)\in(\P^1(\C))^p\times {\hat U_i}$.
Then there is a constant $C'>0$ such that
\begin{equation}
\label{5.17}
\|\sigma(x)\|^2 \leqq C' \frac{|P_i(\lambda\cdot x)|^2}
{\prod_j (1+|u_j|^2)^{d_j}}
\end{equation}
for all $x\in\hat M$, $\lambda\in\Lambda_0$ with $\lambda\cdot x\in
U_i$.
 From (\ref{5.16}) and (\ref{5.17}) it follows that
for $\re z^\rho >0$
\begin{align*}
\|\sigma(f(z))\|^2 &\leqq 
C'C^2\frac{|e^{(\sum_j d_jc_j)z^\rho-c_1z^\rho}|^2}
            {\prod_j(1+|e^{2c_jz^\rho}|)^{d_j}}
\leqq
C'C^2\frac{|e^{(\sum_j d_jc_j)z^\rho -c_1z^\rho}|^2}
            {\prod_j |e^{2c_jd_jz^\rho}|} \\
& = 
C' C^2 |e^{-c_1z^\rho}|^2
= C'C^2e^{-2c_1\re z^\rho}.
\end{align*}
Hence,
$$
\log^+\frac{1}{\|\sigma(f(z))\|} \geqq c_1\re z^\rho + O(1)
$$
for all $z\in\C$ with $\re z^\rho > 0$.
Therefore,
\begin{align*}
m_f(r; \bar D)&=\frac{1}{2\pi}\int_{\{|z|=r\}}
                 \log\frac{1}{\|\sigma(f(z))\|}d\theta\\
        &=\frac{1}{2\pi}\int_{\{|z|=r\}}
                 \log^+\frac{1}{\|\sigma(f(z))\|}d\theta + O(1)\\
        &\geqq\frac{1}{2\pi}\int_{\{|z|=r;\, \re z^\rho>0\}}
                 \log^+\frac{1}{\|\sigma(f(z))\|}d\theta + O(1)\\
        &= \frac{1}{2\pi}\int_{{\{|z|=r\}}}
           c_1\cdot(\re z^\rho)^+ d\theta +O(1) \\
        &= \frac{1}{2\pi}\int_{0}^{2\pi}
           c_1r^\rho \cos^+\rho\theta d\theta +O(1) \\
        &= \frac{c_1}{\pi}r^\rho +O(1).
\end{align*}
On the other hand one deduces easily from (\ref{5.15})
that $T_f(r; D)=O(r^\rho)$.
Hence,
$$
\delta(f; \bar D)=\lowlim_{r\to\infty}
\frac{m_f(r; \bar D)}{T_f(r; \bar D)}>0.
$$
\end{proof}

We will now give an explicit example with $\StD=\{0\}$.

\begin{example}
\label{5.18}{\rm\hskip5pt
Let $A$ be the semi-abelian variety $A=\C^*\times\C^*$,
compactified by $\P^1(\C)\times\P^1(\C)$
with a pair of homogeneous coordinates, $([x_0:x_1],[y_0:y_1])$.
For a pair of natural numbers $(m,n)$ with $m<n$,
let $\bar D$ be the divisor given by
$$
\bar D= \{ ([x_0:x_1],[y_0:y_1]) : y_0^nx_1 +
y_0^{n-m}y_1^mx_0+y_1^nx_0=0 \}.
$$
Set $D=\bar D \cap A$.
Note that $\StD=\{0\}$.
Moreover, $D$ violates condition \ref{4.11},
since $\bar D \ni ([1:0], [1:0])$.
Let $c$ be a positive irrational real number such that
\begin{equation}
\label{5.19}
0<cm<1<cn.
\end{equation}
Let $f:\C\to A$ be the holomorphic curve given by
$$
f: z \mapsto ([1:e^z],[1:e^{cz}]).
$$
Let $\Omega_i$, $i=1,2$, be the Fubini-Study metric forms
of the two factors of $(\P^1(\C))^2$.
Then $c_1(\bar D)=\Omega_1+n\Omega_2$.
By an easy computation one obtains
\begin{equation}
\label{5.20}
T_f(r; c_1(\bar D))  =\frac{1+nc}{\pi} r +O(1).
\end{equation}
Thus, $\rho_f=1$,
and the image $f(\C)$ is Zariski dense in $A$, because $c$ is irrational.

We compute $N(r; f^*D)$ as follows.
Note the following identity for divisors on $\C$:
\begin{equation}
\label{5.21}
f^*D =(e^z+e^{mcz}+e^{ncz})_0.
\end{equation}
We consider a holomorphic curve $g$ in $\P^2(\C)$ with
the homogeneous coordinate system $[w_0: w_1: w_2]$ defined by
$$
g: z \in \C \to [e^z: e^{mcz} : e^{ncz}] \in \P^2(\C).
$$
By computing the Wronskian of $e^z, e^{mcz}$ and $e^{ncz}$
one sees that they are linearly independent over $\C$;
that is $g$ is linearly non-degenerate.
Let $T_g(r)$ be the order function of $g$ with respect to
the Fubini-Study metric form on $\P^2(\C)$.
It follows that
\begin{align}
\label{5.22}
T_g(r) &= \frac{1}{4\pi}\int_{\{|z|=r\}}
\log \left(|e^z|^2+ |e^{mcz}|^2 + |e^{ncz}|^2\right)d\theta +O(1)\\
&= \frac{1}{4\pi}\int_{\{|z|=r\}}
\log\left(1 + |e^{(mc-1)z}|^2 + |e^{(nc-1)z}|^2\right)d\theta+O(1). \nonumber
\end{align}
If $\re z \geqq 0$ (resp.\ $\leqq 0$),
$|e^{(mc-1)z}|\leqq 1$ (resp.\ $\geqq 1$) and
$|e^{(nc-1)z}|\geqq 1$ (resp.\ $\leqq 1$).
Therefore, if $z=re^{i\theta}$ and $\re z \geqq 0$,
\begin{align*}
\log\left(1 + |e^{(mc-1)z}|^2 + |e^{(nc-1)z}|^2\right) &=
2\log^+ |e^{(nc-1)z}|+O(1) \\
&= 2(nc-1)r\cos \theta +O(1).
\end{align*}
If $z=re^{i\theta}$ and $\re z \leqq 0$,
\begin{align*}
\log\left(1 + |e^{(mc-1)z}|^2 + |e^{(nc-1)z}|^2\right) &=
2\log^+ |e^{(mc-1)z}|+O(1) \\
&= 2(mc-1)r\cos \theta +O(1).
\end{align*}
Combining these with (\ref{5.22}), we have
\begin{equation}
\label{5.23}
T_g(r)=\frac{(n-m)c}{\pi} r +O(1).
\end{equation}
We consider the following four lines $H_j$,
$1\leqq j \leqq 4$, of $\P^2(\C)$
in general position:
$$
H_j=\{w_{j-1}=0\} ,\quad 1 \leqq j \leqq 3,\quad  H_4=\{w_0+w_1+w_2=0\}.
$$
Noting that $g$ is linearly non-degenerate and has
a finite order (in fact, $\rho_g=1$), we infer from Cartan's S.M.T. [Ca33] that
\begin{equation}
\label{5.24}
T_g(r)\leqq \sum_{j=1}^4 N_2(r; g^*H_j) +O(\log r).
\end{equation}
Since $N_2(r; g^*H_j)=0$, $1 \leqq j \leqq 3$,
we deduce from (\ref{5.24}), (\ref{5.23}) and (\ref{2.1}) that
$$
N(r; g^*H_4)=\frac{(n-m)c}{\pi} r +O(\log r).
$$
By (\ref{5.21}), $N(r; g^*H_4)=N(r; f^*D)$, and so
\begin{equation}
\label{5.25}
N(r; f^*D)=\frac{(n-m)c}{\pi} r +O(\log r).
\end{equation}
It follows from (\ref{5.20}) and (\ref{5.25}) that
\begin{equation}
\label{5.26}
\delta(f; \bar D) = \frac{1+mc}{1+nc}.
\end{equation}

By elementary calculations one shows that $\ord_z f^*D \geqq 2$ implies
\begin{equation*}
(mc-1)\left( e^{cz}\right)^m + (nc-1) \left( e^{cz}\right)^n=0.
\end{equation*}
Furthermore, $f(z)\in D$ if and only if $e^z+e^{mcz}+e^{ncz}=0$.
Combined, these two relations imply that there is a finite subset
$S\subset\C^2$ such that $\ord_z f^*D \geqq 2$ implies
$(e^z,e^{cz})\in S$.
Since $z\mapsto(e^z,e^{cz})$ is injective, it follows that
$\{z: \ord_z f^*D \geqq 2\}$ is a finite set.
Therefore,
\begin{align}
\label{5.27}
N_1(r,f^*D) &= N(r,f^*D) + O(\log r), \\
\delta_1(f; \bar D) &=\delta(f; \bar D)=\frac{1+mc}{1+nc}.
\nonumber
\end{align}

Let $c'>1$ be an irrational number, and set
$$
c=1/c', \quad m=[c'],\quad n=[c']+1,
$$
where $[c']$ denotes the integral part of $c'$.
Then $m, n$ and $c$ satisfy (\ref{5.19}), and by (\ref{5.26})
$$
\delta(f; \bar D)=
\frac{1+[c']/c'}{1+([c']+1)/c'} \to 1 \qquad (c' \to \infty).
$$
Thus $\delta(f; \bar D)$ ($=\delta_1(f; \bar D)$ by (\ref{5.27}))
takes values arbitrarily close to $1$.
}
\end{example}

\begin{remark}
\label{5.28}{\rm\hskip5pt
In [No96], the first author proved that for $D$ without condition \ref{4.11}
a holomorphic curve $f:\C \to A$, omitting $D$, has no Zariski
dense image, and is contained in a translate of a semi-Abelian
subvariety which has no intersection with $D$.
What was proved in [No96] applied to $f: \C \to A$ with
Zariski dense image yields that there is a positive constant $\kappa$
such that
\begin{equation}
\label{5.29}
\kappa T_f(r; c_1(\bar D))\leqq N_1(r; f^*D) +
S_f(r; c_1(\bar D)),
\end{equation}
provided that $\StD=\{0\}$.
The above $\kappa$ may be, in general, very small because of the
method of the proof.
One needs more detailed properties of $J_k(D)$ to get
the best bound such as in the Main Theorem than
to get (\ref{5.29}); this is the reason why we need the boundary
condition \ref{4.11} for $D$.}
\end{remark}

\begin{remark}
\label{5.30}{\rm\hskip5pt
In [SiY97] Siu and Yeung claimed (\ref{1.2}) for Abelian $A$
of dimension $n$.
The most essential part of their proof was Lemma 2 of [SiY97],
but the claimed assertion does not hold.
We recall the lemma.
\medskip

L{\sc emma} 2 ([SiY97], p.\ 1147). {\it
Let $A$ be an $n$-dimensional Abelian variety and
let $D$ be an ample divisor on $A$.
Let $f: \C \to A$ be a holomorphic curve with Zariski dense image,
and define $W_k$ as in Lemma \ref{5.4}.
Let $\delta \geqq 1$ and $k\geqq n$ be arbitrarily fixed integers.
Then for an arbitrary positive integer $q$
there exists a positive integer $m_0(W_k,\delta,q)$
depending on $W_k, \delta, q$ (and $A$ and $D$) such that
for $m \geqq m_0(W_k,\delta,q)$ there exists an $(L(D))^\delta$-valued
holomorphic $k$-jet differential $\omega$ on $A$ of weight $m$ whose
restriction to $A\times W_k$ is not identically zero and which
vanishes along $J_k(D)\cap (A\times W_k)$ to order at least $q$.
In particular, from the definition of $W_k$ one knows that
$\omega$ is not identically zero on $J_k(f)$.
}

In the proof of (\ref{1.2}) they applied this lemma, taking
$\delta=1$ and $k\geqq n$ fixed, and increasing $q \to \infty$.

Since $A\times W_k$ is the Zariski closure of
the transcendental holomorphic mapping
$J_k(f):\C \to J_k(A)\cong A\times \C^{nk}$,
the variety $W_k$ must be allowed to be quite arbitrary.
In fact, in the proof of the above Lemma 2, [SiY97],
the fact that $A\times W_k$ was defined to be the Zariski closure
of the $k$-jet lifting of $f:\C \to A$ with Zariski dense image
was not used at all except for the last statement, and hence Lemma 2
(except for the last statement) should be true for
arbitrary non-empty subvariety $W_k \subset \C^{nk}$
and ample $D$ on $A$,
if the proof were correct.
This is a very different point from our proof (cf. the proof of
Lemma \ref{5.4}).
But, we then deduce some contradictory conclusions as follows.

(a)  We take an ample divisor $D$ on $A$ such that it contains
a translate of a non-trivial Abelian subvariety $A'$
(cf.\ [\NO], Example (6.4.13) for such an example).
Let $g:\C \to A'$ be a one-parameter subgroup with
Zariski dense image.
We regard $g$ as a holomorphic curve into $A$, and set
$f(z)=g(z)+a$ with $a \in A \setminus D$.
Then $f$ is a holomorphic curve such that $f(\C) \not\subset D$,
and  $W_k$ consists of only one point for every $k \geqq 1$.
We obtain $A\times W_k\cong A$.
Through this isomorphism, we have that
$$
A' \subset J_k(D) \cap (A\times W_k) \subsetneqq A.
$$
Let $\mathcal{I}_k=\mathcal{I}(J_k(D) \cap (A\times W_k))$
denote the ideal sheaf of
$J_k(D) \cap (A\times W_k)$ ($\subset A$).
Note that any jet differential of any weight $m$ restricted to
$A\times W_k$ is reduced to a jet differential of weight $0$,
for $W_k$ consists of one point.
Then Lemma 2 should imply that for all $q \geqq 1$
$$
H^0(A, \mathcal{O}((L(D))^\delta) \tensor \mathcal{I}_k^q)\not=\{0\},
$$
where $\mathcal{O}((L(D))^\delta)$ denotes the sheaf of germs of
holomorphic sections of $(L(D))^\delta$.
Since $A' \subset J_k(D) \cap (A\times W_k)$,
the ideal sheaf $\mathcal{I}=\mathcal{I}(A')$ of $A'$
contains $\mathcal{I}_k$.
Therefore,
$$
H^0(A, \mathcal{O}((L(D))^\delta) \tensor \mathcal{I}^q)\not=\{0\},
\qquad \forall q \geqq 1,
$$
and hence the infinite dimensionality of
$H^0(A, \mathcal{O}((L(D))^\delta))$ would follow,
where $\delta$ had been fixed.
This is clearly absurd.
This observation implies that the Zariski denseness
of the image $f(\C)$ in $A$ must be used essentially.

(b)  We also observe that Lemma 2 is not valid even for
$f:\C \to A$ with Zariski dense image,
and moreover that $k$ cannot be fixed as stated in Lemma 2.
Let $k\geqq n$ be any fixed.
Let $f:\C \to A$ be a one-parameter subgroup with Zariski dense image.
Let $D$ be an ample divisor on $A$ containing the zero $0 \in A$
such that $f(\C)$ is tangent highly enough to $D$ at $0$ so that
$J_kf(0) \in J_k(D)$, but $f(\C) \not\subset D$.
Let $\mathfrak{m}_0$ be the maximal ideal sheaf of the structure
sheaf $\mathcal{O}_A$ at $0$.
Since $W_k$ consists of only one point,
$A\times W_k \cong A$, and through this isomorphism
$0\in J_k(D)\cap (A\times W_k)$.
Therefore we have that $\mathcal{I}_k=\mathcal{I}(J_k(D)\cap (A\times W_k))
\subset \mathfrak{m}_0$.
As in (a), Lemma 2 should imply that
$$
H^0(A, \mathcal{O}((L(D))^\delta) \tensor \mathfrak{m}_0^q)\supset
H^0(A, \mathcal{O}((L(D))^\delta) \tensor \mathcal{I}_k^q)\not=\{0\},
\quad \forall q \geqq 1.
$$
Thus we would obtain that $\dim H^0(A, \mathcal{O}((L(D))^\delta))=\infty$;
this is a contradiction.

(c)  The reason of the contradictions observed in
(a) and (b) with respect to the above Lemma 2 comes from
the use of the semi-continuity theorem for a non-flat
family of coherent ideal sheaves.
They used a deformation technique of the given
ample divisor $D$.
That is, taking a generic small deformation family
$D(t), t \in \Delta(1)$ on $A$ with $D(0)=D$,
they considered the family of ideals,
$\left\{\big(\mathcal{I}(J_k(D(t))\cap(A\times W_k))
\big)^q\right\}_{t \in\Delta(1)}$
for $k\geqq n$ and $q\geqq 1$.
More precisely, they worked on the compactification
$\overline{J_k(A)}=A \times \P^{nk}(\C)$
of $J_k(A)\cong A \times \C^{nk}$.
Let $\overline{J_k(D(t))}$ (resp.\ $\bar W_k$) denote
the closure of $J_k(D(t))$ (resp.\ $W_k$) in
$\overline{J_k(A)}$ (resp.\ $\P^{nk}(\C)$).
To apply the semi-continuity theorem of the dimension
of cohomology groups, one needs the flatness of the family,
$\{\mathcal{I}(\overline{J_k(D(t))}\cap(A\times \bar W_k))\}_{t\in \Delta(1)}$.
In general, the constructed family
$\{\mathcal{I}(\overline{J_k(D(t))}\cap(A\times \bar W_k))\}_{t\in \Delta(1)}$
may not be flat,
since there may be a ``jump'' of the supports of those ideals.
This fact tells us the difficulty to apply the
deformation technique to obtain the second main theorem in general.
Because of its own interest we give such an example in what follows.

Let $E$ be an elliptic curve defined by the square lattice,
$\Z+i\Z$, and set $A=E \times E$.
Let $(x,y)$ be a local flat coordinate system of $A$,
and define a holomorphic curve $f:\C \to A$ by
$$
f: z \in \C \to (z, \alpha z) \in A,
$$
where $\alpha$ is an irrational number.
Then the image $f(\C)$ is Zariski dense in $A$.
The 2-jet lifting of $f$ is given by
$$
J_2(f)(z)=((z, \alpha z), (1, \alpha), (0,0))\in
J_2(A) \cong A \times \C^2 \times \C^2.
$$
Thus,
$$
W_2=((1, \alpha), (0,0)).
$$
Let $L$ be a sufficiently ample line bundle over $A$ such that
$L$ carries a global holomorphic section $\sigma(x,y)$
whose germ at $(0,0)$ is written as
$$
- y^2 +x^3+\alpha^2 x^2 + x^4G(x,y).
$$
Let $D$ be the divisor defined by the zero locus of $\sigma$.
By an easy computation one sees
$$
J_2(D)_{(0,0)}=\{(0,0)\} \times  \C(1, \pm \alpha) \times \C^2.
$$
Therefore we have
$$
J_2(D) \cap (A \times W_2) \ni ((0,0), (1, \alpha), (0,0)).
$$
For small $t \in \C$ we consider a generic deformation
$D(t)$ defined by
$$
- y^2 +x^3+\alpha^2 x^2 + x^4G(x,y)+tH(x,y)=0.
$$
We look for a point $(x_0,y_0) \in D(t)$ near $(0,0)$
with $t \not= 0$
such that
$(x_0, y_0)\times (1,\alpha) \times (0,0) \in J_2(D(t))\cap(A\times W_2)$.
First we have
\begin{equation}
\label{5.31}
- y_0^2 +x_0^3+\alpha^2 x_2^2 + x_0^4G(x_0,y_0)+tH(x_0,y_0)=0.
\end{equation}
Set
\begin{align*}
\phi(z)= &- (y_0+\alpha z)^2 +(x_0+z)^3+\alpha^2 (x_0+z)^2
 + (x_0+z)^4G(x_0+z,y_0+\alpha z)\\
&+tH(x_0+z,y_0+\alpha z).
\end{align*}
Then one gets
\begin{align*}
\phi'(z)=&- 2\alpha(y_0+\alpha z) + 3(x_0+z)^2+ 2\alpha^2 (x_0+z)
 + (x_0+z)^3G_1(x_0+z,y_0+\alpha z)\\
&+
tH_x(x_0+z,y_0+\alpha z)+ t\alpha H_y(x_0+z,y_0+\alpha z),
\end{align*}
where $G_1$ is a naturally defined holomorphic function.
Hence,
\begin{align}
\label{5.32}
\phi'(0)=&- 2\alpha y_0 + 3x_0^2+ 2\alpha^2 x_0
 + x_0^3G_1(x_0,y_0)\\
&+tH_x(x_0,y_0)+ t\alpha H_y(x_0,y_0)=0.\nonumber
\end{align}
Taking the second derivative, we have
\begin{align*}
\phi''(z)= &6(x_0+z)
 + (x_0+z)^2G_2(x_0+z,y_0+\alpha z)\\
&+tH_{xx}(x_0+z,y_0+\alpha z)
+ 2t\alpha H_{xy}(x_0+z,y_0+\alpha z)
+ t\alpha^2 H_{yy}(x_0+z,y_0+\alpha z),
\end{align*}
and so
\begin{align}
\label{5.33}
\phi''(0)=&(6+ x_0 G_2(x_0,y_0)) x_0\\
&+(H_{xx}(x_0,y_0)+ 2\alpha H_{xy}(x_0,y_0)+\alpha H_{yy}(x_0,y_0))t=0.
\nonumber
\end{align}
We may assume that for
$(x_0,y_0)$ close to $(0,0)$
$$
6+ x_0 G_2(x_0,y_0)\not=0.
$$
Thus we may write
$$
x_0=t\psi(t,y_0).
$$
We substitute this to (\ref{5.32}), and get
\begin{align*}
&- 2\alpha y_0 + 3t^2\psi^2(t,y_0)+ 2\alpha^2 t\psi(t,y_0)
 + t^3\psi^3(t,y_0)G_1(t\psi(t,y_0),y_0)\\
&+tH_x(t\psi(t,y_0),y_0)+ t\alpha H_y(t\psi(t,y_0),y_0)=0.
\end{align*}
Therefore we have
$$
y_0=t\lambda (t),\quad x_0=t\psi(t,t \lambda(t))=t\mu(t).
$$
Then we substitute these to (\ref{5.31}), and obtain
$$
- t^2\lambda^2 (t) +t^3\mu^3 (t)
+\alpha^2 t^2\mu^2(t) + t^4\mu^4(t)G(t\mu(t),t\lambda (t))
+tH(t\mu(t),t\lambda (t))=0.
$$
Since $t\not=0$, we have
\begin{equation}
\label{5.34}
- t\lambda^2 (t) +t^2\mu^3 (t)
+\alpha^2 t\mu^2(t) + t^3\mu^4(t)G(t\mu(t),t\lambda (t))
+ H(t\mu(t),t\lambda (t))=0.
\end{equation}
We assume a generic condition, $H(0,0) \not=0$;
equation (\ref{5.34}) is not trivial.
Thus, $t$ satisfying (\ref{5.34}) is isolated, and cannot approach $0$.

It follows that there is a neighborhood $U \subset J_2(A)$ of
$((0,0), (1, \alpha), (0,0)) \in J_2(D(0)) \cap (A \times W_2)$
such that for every small $t \not= 0$,
$$
J_2(D(t)) \cap (A \times W_2) \cap U= \emptyset.
$$
Therefore the ideal family
$\{\mathcal{I}(\overline{J_k(D(t))}\cap(A\times \bar W_k))\}_{t\in \Delta(1)}$
is not flat.
}
\end{remark}

\begin{remark}
\label{5.35}{\rm\hskip5pt
It is an interesting problem to see if
the truncation level
$k_0$ of the counting function $N_{k_0}(r; f^*D)$
in the Main Theorem can be taken as a function only in $\dim A$.
By the above proof, it would be sufficient to find a natural number
$k$ such that
$\pi_2(J_k(\bar D; \log \bar A \cap \partial A))\cap W_k \not= W_k$.
Note that $\dim \pi_2(J_k(\bar D; \log \bar A))\leqq
\dim J_k(\bar D; \log \bar A)=(n-1)(k+1)$.
Thus, if $\dim W_k > (n-1)(k+1)$ we may set $k_0=k$.
For example, if $J_n(f)(\mathbf{C})$ is Zariski dense in $J_n(A)$,
then $\dim W_n=n^2$.
Since $\dim \pi_2(J_n(\bar D; \log \partial A))=n^2-1$,
we may set $k_0=n$.
}
\end{remark}

\section{Applications.}

Let the notation be as in the previous section.
Here we assume that $A$ is an Abelian variety and
$D$ is reduced and hyperbolic;
in this special case, $D$ is hyperbolic if and only if $D$ contains
no translate of a one-parameter subgroup of $A$.
Cf.\ [\NO], [La87] and [Ko98] for the theory of hyperbolic complex
spaces.

\begin{theorem}
\label{6.1}
Let $D \subset A$ be hyperbolic and $d_0$ be the highest order
of tangency of $D$ with translates of one-parameter
subgroups.
Let $\pi: X \to A$ be a finite covering space such that its ramification
locus contains $D$ and the ramification order over $D$ is
greater than $d_0+1$.
Then $X$ is hyperbolic.
\end{theorem}

{\it Proof.}
By Brody's theorem (cf., e.g., [\NO], Theorem (1.7.3))
it suffices to show that
there is no non-constant holomorphic curve
$g:\mathbf{C} \to X$ such that
the length $\|g'(z)\|$ of the derivative $g'(z)$ of $g(z)$
with respect to an arbitrarily fixed Finsler metric
on $X$ is bounded.
Set $f(z)=\pi(g(z))$.
Then the length $\|f'(z)\|$ with respect to the
flat metric is bounded, too, and hence
$f'(z)$ is constant.
Thus, $f(z)$ is a translate of a one-parameter subgroup.
By definition we may take $k_0=d_0+1$ in (\ref{5.7}).
Take $d\ (> d_0+1)$ so that $X$ ramifies over $D$ with order at least $d$.
Then we have that
$N_1(r; f^*D) \leqq \frac{1}{d}N(r; f^*D)$.
Hence it follows from the Main Theorem that
\begin{align*}
T_f(r; L(D)) &= N_{d_0+1}(r, f^*D) + O(\log r)
\leqq (d_0+1) N_1(r; f^*D) + O(\log r)\\
&\leqq \frac{d_0+1}{d}N(r; f^*D) + O(\log r)
\leqq \frac{d_0+1}{d}T_f(r; L(D)) + O(\log r).
\end{align*}
Since $T_f(r; L(D)) \geqq c_0 r^2$ with a constant $c_0>0$,
$d \leqq d_0+1$;
this is a contradiction.
{\it Q.E.D.}
\medskip

{\it Remark.}
In the special case of $\dim X=\dim A=2$,
C.G. Grant [Gr86] proved that if $X$ is of general type
and $X\to A$ is a finite (ramified) covering space,
then $X$ is hyperbolic.
When $\dim X=\dim A=2$, $D$ is an algebraic curve,
and hence the situation is much simpler than
the higher dimensional case.

\begin{theorem}
\label{6.2}
Let $f:\mathbf{C} \to A$ be a 1-parameter analytic subgroup
in $A$ with $a=f'(0)$.
Let $D$ be an effective divisor on $A$ with the Riemann form
$H(\cdot,\cdot)$.
Then we have
$$
N(r; f^*D)=H(a,a)\pi r^2 +O(\log r).
$$
\end{theorem}

{\it Proof.}
Note that the first Chern class $c_1 (L(D))$ is represented by
$i \partial \bar\partial H(w,w)$.
It follows from (\ref{2.1}) and Lemma \ref{5.1} that
\begin{align*}
N(r; f^*D) &=T_f(r; L(D))+O(\log r)\\
&=\int_0^r \frac{dt}{t}\int_{\Delta(t)}
i H(a,a) dz \wedge d\bar z+O(\log r)\\
&=H(a,a) \pi r^2 +O(\log r).
\end{align*}
{\it Q.E.D.}
\medskip

\begin{remark}
\label{6.3}{\rm\hskip5pt
In the case where $f(\C)$ is Zariski dense in $A$,
Ax ([Ax72]) proved the following estimate,
$$
0 < \lowlim_{r\to\infty} \frac{n(r, f^*D)}{r^2} \leqq
\upplim_{r\to\infty}\frac{n(r, f^*D)}{r^2} <\infty,
$$
which is equivalent to
$$
0 < \lowlim_{r\to\infty} \frac{N(r, f^*D)}{r^2} \leqq
\upplim_{r\to\infty}\frac{N(r, f^*D)}{r^2} <\infty.
$$
}
\end{remark}

\section*{\centerline{\normalsize\it References}}

\parindent35pt
\begin{itemize}
\setlength{\itemsep}{-3pt}
\item[{[Ah41]}]
L. V. Ahlfors,
The theory of meromorphic curves,
Acta Soc.\ Sci.\ Fennicae, Nova Ser.\ A.
{\bf 3}
(1941),
3-31.
\item[{[AN91]}]
Y. Aihara and J. Noguchi,
Value distribution of meromorphic mappings into compactified locally symmetric spaces,
Kodai Math.\ J.
{\bf 14}
(1991),
320-334.
\item[{[Ax72]}]
J. Ax,
Some topics in differential algebraic geometry II,
Amer.\ J. Math.\
{\bf 94}
(1972),
1205-1213.
\item[{[Bl26]}]
Bloch, A.:
Sur les syst\`emes de fonctions uniformes satisfaisant \`a l'\'equation d'une vari\'et\'e alg\'ebrique dont l'irr\'egularit\'e d\'epasse la dimension.
J. Math.\ Pures Appl.\
{\bf 5}
(1926),
19-66.
\item[{[CG72]}]
J. Carlson and P. Griffiths,
A defect relation for equidimensional holomorphic
mappings between algebraic varieties,
Ann.\ Math.\ {\bf 95} (1972), 557-584.
\item[{[Ca33]}]
H. Cartan,
Sur les z\'eros des combinaisons lin\'eaires de $p$ fonctions holomorphes
donn\'ees, Mathematica \textbf{7} (1933), 5-31.
\item[{[DL97]}]
G. Dethloff and S.S.Y. Lu,
Logarithmic projective jet bundles and applications,
preprint, 1997.
\item[{[ES92]}]
A.E. Eremenko and M.L. Sodin,
The value distribution of meromorphic functions and meromorphic curves
from the view point of potential theory,
St.\ Petersburg Math.\ J.\ {\bf 3} (1992) No.\ 1, 109-136.
\item[{[Gr86]}]
C.G. Grant, Entire holomorphic curves in surfaces,
Duke Math. J. {\bf 53} (1986), 345-358.
\item[{[GK73]}]
P. Griffiths and J. King,
Nevanlinna theory and holomorphic mappings between algebraic varieties,
Acta Math.\
{\bf 130} (1973), 145-220.
\item[{[Ha64]}]
W.K. Hayman, Meromorphic Functions, Oxford Math. Monographs,
Oxford University Press, London, 1964.
\item[{[Kr98]}]
R. Kobayashi, Holomorphic curves in Abelian varieties:
The second main theorem and applications, preprint, 1998.
\item[{[Ko98]}]
S. Kobayashi, Hyperbolic Complex Spaces, Grundlehren der mathematischen
Wissenschaften {\bf 318}, Springer-Verlag, Berlin-Heidelberg, 1998.
\item[{[La87]}]
S. Lang,
Introduction to Complex Hyperbolic Spaces,
Springer-Verlag, New York-Berlin-Heidelberg, 1987.
\item[{[Mc96]}]
M. McQuillan,
A dynamical counterpart to Faltings' ``Diophantine approximation on
Abelian varieties'',
I.H.E.S. preprint, 1996.
\item[{[No77]}]
J. Noguchi,
Holomorphic curves in algebraic varieties,
Hiroshima Math.\ J.
{\bf 7}
(1977),
833-853.
\item[{[No81]}]
J. Noguchi,
Lemma on logarithmic derivatives and holomorphic curves in algebraic varieties,
Nagoya Math.\ J.
{\bf 83}
(1981),
213-233.
\item[{[No86]}]
J. Noguchi,
Logarithmic jet spaces and extensions of de Franchis' theorem,
Contributions to Several Complex Variables,
pp.\ 227-249,
Aspects Math.\ No.\
{\bf 9},
Vieweg,
Braunschweig,
1986.
\item[{[No96]}]
J. Noguchi,
On Nevanlinna's second main theorem,
Geometric Complex Analysis,
Proc.\ the Third International Research Institute, Math.\ Soc.\ Japan,
Hayama, 1995, pp.\ 489-503, World Scientific, Singapore, 1996.
\item[{[No98]}]
J. Noguchi,
On holomorphic curves in semi-Abelian varieties,
Math.\ Z. {\bf 228} (1998), 713-721.
\item[{[NW99]}]
J. Noguchi and J. Winkelmann,
Holomorphic Curves and Integral Points off Divisors, preprint, 1999.
\item[{[\NO]}]
J. Noguchi and T. Ochiai,
Geometric Function Theory in Several Complex Variables,
Japanese edition, Iwanami, Tokyo, 1984;
English Translation, Transl.\ Math.\ Mono.\ {\bf 80},
Amer.\ Math.\ Soc., Providence, Rhode Island,
1990.
\item[{[Si87]}]
Y.-T. Siu, Defect relations for holomorphic maps between spaces of
different dimensions,
Duke Math.\ J.\  {\bf 55} (1987), 213-251.
\item[{[SiY96]}]
Y.-T. Siu and S.-K. Yeung,
A generalized Bloch's theorem and the hyperbolicity of the complement
of an ample divisor in an Abelian variety,
Math.\ Ann.\ {\bf 306} (1996), 743-758.
\item[{[SiY97]}]
Y.-T. Siu and S.-K. Yeung,
Defects for ample divisors of Abelian varieties, Schwarz lemma, and
hyperbolic hypersurfaces of low degrees,
Amer.\ J.\ Math.\ {\bf 119} (1997), 1139-1172.
\item[{[St53/54]}]
W. Stoll,
Die beiden Haupts\"atze der Wertverteilungstheorie bei
Funktionen mehrerer komplexer Ver\"anderlichen (I); (II),
Acta Math.\ \textbf{90} (1953), 1-115;
ibid,
Acta Math.\ \textbf{92} (1954), 55-169.
\item[{[We58]}]
Weil, A.:
Introduction \`a l'\'etude des vari\'et\'es k\"ahleriennes,
Hermann, Paris,
1958.
\end{itemize}

\bigskip
\baselineskip=12pt
\rightline{Graduate School of Mathematical Sciences}
\rightline{University of Tokyo}
\rightline{Komaba, Meguro,Tokyo 153-8914}
\rightline{e-mail: noguchi@ms.u-tokyo.ac.jp}
\bigskip

\rightline{Mathematisches Institut}
\rightline{Rheinsprung 21}
\rightline{CH--4053 Basel}
\rightline{Switzerland}
\rightline{e-mail: jwinkel@member.ams.org} 
\rightline{URL: http://www.cplx.ruhr-uni-bochum.de/$\tilde{\ }$jw/index-e.html}
\bigskip

\rightline{Research Institute for Mathematical Sciences}
\rightline{Kyoto University}
\rightline{Oiwake-cho, Sakyoku, Kyoto 606-8502}
\rightline{e-mail: ya@kurims.kyoto-u.ac.jp}

\end{document}